\documentclass{amsart}
\usepackage{verbatim}
\usepackage{amsthm,amsmath,amsfonts}
\usepackage{mathrsfs}
\usepackage[colorlinks=true,linktocpage=true]{hyperref}
\usepackage{tocvsec2}
\usepackage{accents}
\usepackage{hyperref}
\usepackage{upgreek}
\usepackage{stmaryrd}
\usepackage[usenames,dvipsnames]{xcolor}
\hypersetup{urlcolor=cyan, citecolor=cyan, linkcolor=red}

\usepackage{framed}
\usepackage{bbm,dsfont}
\usepackage{color}
\numberwithin{equation}{section}
\numberwithin{table}{section}
\numberwithin{figure}{section}
\font\tenscrpt=eusm10
\font\sevenscrpt=eusm10 scaled 700
\font\fivescrpt=eusm10 scaled 500
\newfam\eusmfam
\textfont\eusmfam=\tenscrpt
\scriptfont\eusmfam=\sevenscrpt
\scriptscriptfont\eusmfam=\fivescrpt

\newtheorem{thm}{Theorem}[section]
\newtheorem{cor}{Corollary}[section]
\newtheorem{lem}{Lemma}[section]
\newtheorem{prop}{Proposition}[section]

\theoremstyle{definition}
\newtheorem{defn}{Definition}[section]

\newtheorem{rem}{Remark}[section]
\newtheorem{notn}{Notation}[section]

\newtheorem{ack}{Acknowledgment\hspace{-1.5mm}}

\newcommand{\thmref}[1]{Theorem~\ref{#1}}
\newcommand{\secref}[1]{Section~\ref{#1}}

\newcommand{\appref}[1]{Appendix~\ref{#1}}
\newcommand{\lemref}[1]{Lemma~\ref{#1}}

\newcommand{\remref}[1]{Remark~\ref{#1}}

\newcommand{\figref}[1]{Figure~\ref{#1}}
\def\qed{\quad\vcenter{\hrule\hbox{\vrule height.6em\kern.6em\vrule}\hrule}}
\newenvironment{pf}{{\bigskip\textit{\newline Proof.}\quad}}{$\qed$\bigskip\newline}
\newenvironment{pf*}[1]{{\bigskip\textit{\newline#1.}\quad}}{$\qed$\bigskip\newline}
\def\ds{\displaystyle}
\def\ts{\textstyle}

\def\Cmpl{\mathbb C}

\def\psxy{{K}^{\text{\tiny{\sc{BM}}}^d}_{s;x,y}}

\def\K{{\mathbb K}}

\def\KBtxy{{\K}^{\text{\tiny{\sc{BTBM}}}^d}_{t;x,y}}

\def\KBtx{{\K}^{\text{\tiny{\sc{BTBM}}}^d}_{t;x}}

\def\KBtdot{{\K}^{\text{\tiny{\sc{BTBM}}}^d}_{t;\cdot}}

\def\KKStxy{{\K}^{\text{\tiny{\sc{nlksa}}}^d}_{t;x,y}}

\def\KKStsxy{{\K}^{\text{\tiny{\sc{nlksa}}}^d}_{t-s;x,y}}

\def\KKtsSy{{\K}^{(\beta,d)}_{t-s;\cdot,y}}

\def\psxy{{K}^{\text{\tiny{\sc{BM}}}^d}_{s;x,y}}

\def\ptsz{{K}^{\text{\tiny{\sc{BM}}}}_{t;s}}

\def\KIhalftzxsqrt2{K^{\Lambda_{\frac12}}_{t;0,\sqrt2 x}}
\def\KabBmtzxsqrt2{K^{{|B|}}_{t;0,\frac{x}{\sqrt{2}}}}

\def\K{{\mathbb K}}
\def\sKKS{{\mathscr K}^{\text{\tiny{\sc{LKS}}}^d}}

\def\KKtsSy{{\K}^{\text{\tiny{\sc{LKS}}}^d}_{t-s;\cdot,y}}

\def\KKStdot{{\K}^{\text{\tiny{\sc{LKS}}}^d}_{t;\cdot}}

\def\KKStxy{{\K}^{\text{\tiny{\sc{LKS}}}^d}_{t;x,y}}

\def\KKStsxy{{\K}^{\text{\tiny{\sc{LKS}}}^d}_{t-s;x,y}}
\def\KKStrxy{{\K}^{\text{\tiny{\sc{LKS}}}^d}_{t-r;x,y}}

\def\KKSthrdot{{\K}^{\text{\tiny{\sc{LKS}}}^d}_{\theta-r;\cdot}}
\def\KKSthrdoty{{\K}^{\text{\tiny{\sc{LKS}}}^d}_{\theta-r;\cdot,y}}
\def\KKStrdot{{\K}^{\text{\tiny{\sc{LKS}}}^d}_{t-r;\cdot}}

\def\KKStrdoty{{\K}^{\text{\tiny{\sc{LKS}}}^d}_{t-r;\cdot,y}}
\def\KKStrxdot{{\K}^{\text{\tiny{\sc{LKS}}}^d}_{t-r;x,\cdot}}
\def\KKSsrdot{{\K}^{\text{\tiny{\sc{LKS}}}^d}_{s-r;\cdot}}
\def\KKSsrdoty{{\K}^{\text{\tiny{\sc{LKS}}}^d}_{s-r;\cdot,y}}

\def\KKStx{{\K}^{\text{\tiny{\sc{LKS}}}^d}_{t;x}}
\def\KKStzx{{\K}^{\text{\tiny{\sc{LKS}}}^d}_{t;0,x}}

\def\FKKStxi{\hat{\K}^{\text{\tiny{\sc{LKS}}}^d}_{t;\xi}}

\def\i{\mathbf i}

\def\pa{\partial}

\def\Rpop{\mathring{\R}_{+}}
\def\L{\mathbb{L}}
\def\e{\mathrm{e}}

\def\unx{u_0(x)}

\def\uny{u_0(y)}
\def\un{u_0}

\def\N{{\mathbb N}}

\def\Rd{{\mathbb R}^{d}}

\def\Rpop{\mathring{\R}_{+}}
\def\R{\mathbb R}

\def\S{\mathbb S}

\def\X{\mathbb X}
\def\Rp{{\R}_+}

\def\sF{{\mathscr F}}

\def\eqdef{:=}

\def\lap{\Delta}

\def\df#1#2{\ds{\frac{#1}{#2}}}
\def\tf#1#2{\ts{\frac{#1}{#2}}}

\def\lbl#1{\label{#1}}

\def\intrd{\int_{\Rd}}
\def\intrdzt{\int_{\Rd}\int_0^t}

\def\pa{\partial}

\def\lang{\left<}
\def\rang{\right>}
\def\lab{\left|}
\def\rab{\right|}
\def\lpa{\left(}
\def\rpa{\right)}

\def\lbk{\left[}
\def\rbk{\right]}
\def\lbr{\left\{}
\def\rbr{\right\}}
\def\back{\begin{ack}}
\def\eack{\end{ack}}
\def\bdf{\begin{defn}}
\def\edf{\end{defn}}
\def\bcr{\begin{cor}}
\def\ecr{\end{cor}}
\def\bnt{\begin{notn}}
\def\ent{\end{notn}}
\def\brm{\begin{rem}}
\def\erm{\end{rem}}
\def\blm{\begin{lem}}
\def\elm{\end{lem}}
\def\bpf{\begin{pf}}
\def\bpfs{\begin{pf*}}
\def\epf{\end{pf}}
\def\epfs{\end{pf*}}
\def\bfg{\begin{figure}}
\def\efg{\end{figure}}

\def\beq{\begin{equation}}
\def\beqs{\begin{equation*}}
\def\eeq{\end{equation}}
\def\eeqs{\end{equation*}}
\def\beqar{\begin{eqnarray}}
\def\beqars{\begin{eqnarray*}}
\def\eeqar{\end{eqnarray}}
\def\eeqars{\end{eqnarray*}}
\def\bsp{\begin{split}}
\def\esplit{\end{split}}
\def\bc{\begin{cases}}
\def\ec{\end{cases}}
\def\bt{\begin{tabular}}
\def\et{\end{tabular}}

\def\bthm{\begin{thm}}
\def\ethm{\end{thm}}
\def\bpr{\begin{prop}}
\def\epr{\end{prop}}
\def\bfr{\begin{framed}}
\def\efr{\end{framed}}
\def\bfg{\begin{figure}}
\def\efg{\end{figure}}

\def\bsh{\begin{shaded}}
\def\esh{\end{shaded}}
\def\bcm{\iffalse}
\def\babs{\begin{abstract}}
\def\eabs{\end{abstract}}
\def\bit{\begin{itemize}}
\def\eit{\end{itemize}}
\def\ben{\begin{enumerate}}

\def\rencomrom{\renewcommand{\labelenumi}{(\roman{enumi})}}
\def\een{\end{enumerate}}
\def\bdsc{\begin{description}}
\def\edsc{\end{description}}
\def\babs{\begin{abstract}}
\def\eabs{\end{abstract}}

\def\nim{\wedge}

\def\sKKS{{\mathscr K}^{\text{\tiny{\sc{LKS}}}^d}}

\def\sKKSl{{\mathscr K}^{\text{\tiny{\sc{LKS}}}^d}_{\ell}}

\def\sO{\mathscr O}

\def\ig{\iffalse}

\def\e{\mathrm{e}}

\def\i{\mathbf i}
\def\sI{\mathscr I}
\def\sJ{\mathscr J}

\def\pa{\partial}

\def\bt{b^{(2)}}

\def\l{\ell}
\def\B{\mathds B}
\def\C{\mathds C}

\def\Lp{{\mathds L}^{p}}
\def\Lq{{\mathds L}^{q}}

\def\Ltwo{{\mathds L}^{2}}
\def\Lt{{\mathds L}^{2}}
\def\Ltp{{\mathds L}^{2p}}

\def\gamo{\gamma_{1}}
\def\gamt{\gamma_{2}}
\def\Lgam{{\mathds L}^{\gamma}}
\def\Lgamo{{\mathds L}^{\gamma_{1}}}
\def\Lgamt{{\mathds L}^{\gamma_{2}}}

\def\MsB{{\mathds M}_{\mathscr B}}

\def\U{\mathds U}

\def\ugatmr{\upgamma_{t-r}}

\def\Xx{X^x}
\def\Xix{X^{-ix}}
\def\Uk{U^{(k)}}
\def\Uo{U^{(1)}}
\def\uno{\un^{(1)}}
\def\Ut{U^{(2)}}
\def\unt{\un^{(2)}}
\def\unk{\un^{(k)}}

\usepackage{graphicx}
\usepackage{scalerel, pict2e}
\title[Kuramoto-Sivashinsky Burgers PDE until the $6p${\scriptsize-th} dimension]{Solving the Kuramoto-Sivashinsky-Burgers equation until the $6p${\scriptsize-th} dimension: the Brownian-time paradigm}
  
\author{Hassan Allouba}
\address{Department of Mathematical Sciences, Kent State University, Kent,
Ohio 44242}
\email{allouba@math.kent.edu}

\begin{document}
\begin{abstract}
We use our earlier Brownian-time framework to formulate and establish  global uniqueness and local-in-time existence of the Burgers incarnation of the Kuramoto-Sivashinsky PDE on $\Rp\times\Rd$, in the class  of time-continuous $\Ltp$-valued solutions, $p\ge1$,  for every $d<6p$.  We assume neither space compactness, nor spatial coordinates dependence, nor smallness of initial data.  The surprising discovery of the $6p$-th dimension bound, even for local solutions, is revealed by our approach and the Brownian-time  kernel---the Brownian average of an angled $d$-dimensional Schr\"odinger propagator---at its heart.  We use this kernel to give a systematic approach, for all dimensions simultaneously, including a novel formulation---even in the well-known $d=1$ case---of the KS equation.  This yields the estimates leading to this article's conclusions.  We achieve the stated results by fusing some of our earlier Brownian-time stochastic processes constructions and ideas---encoded in the aforementioned kernel---with analytic ones, including complex and harmonic analysis; by employing suitable $N$-ball approximations together with fixed point theory; and by an adaptation of the stochastic analytic stopping-time technique to our deterministic setting.  Using a separate strategy, that is also built on our Brownian-time paradigm, we treat the global wellposedness of the multidimensional KS equation in a followup upcoming article.  This work also serves as a template for another forthcoming article in which we prove similar results for the time-fractional Burgers equation in multidimensional space.
\end{abstract}
\maketitle
\tableofcontents
\section{Introduction and statement of results}
\subsection{The Kuramoto-Sivashinsky PDE and its Burgers incarnation} 
The Kuramoto Sivashinsky (KS) equation, which has garnered the interest of too many experts to cite here (see e.g., Temam's classic \cite{T} and the references therein and \cite{BenKukaRusZian14,ChesFoi,ZGIK03,Kuk92}), is used for models related to turbulence in chemistry and combustion.  In space dimension $d=1$, it was introduced by Kuramoto \cite{Kur1} for the study of phase turbulence in the Belousov-Zhabotinsky reactions.  Sivashinsky \cite{Siv1,Siv2}  extended Kuramoto's equation to $d\ge2$ to investigate the propagation of a flame front in the case of mild combustion. In its standard eikonal form, the KS PDE  is a scalar equation given by\footnote{\lbl{fn:not1} For aesthetics,  typesetting convenience, and/or presentation flow reasons, we alternate freely between $\pa/\pa{x_{k}}$ and $\pa_{k}$---for the partial derivative in the $k$-th spatial variable, $k\in\{1,\ldots,d\}$---and between $\pa_{t}$ and $\pa/\pa t$ for  the partial derivative in the time variable $t$.  The mixed derivative is also denoted by $\pa^{2}_{t,k}$.} 
\beq\lbl{eq:ksst}
\bsp 
&\displaystyle{\partial_{t}\tilde U}=-\tfrac18\lpa\lap^{2}+4\lap\rpa \tilde U -\tfrac12\sum_{k=1}^{d}\lpa{\pa_{k}\tilde U}\rpa^{2}
,\ (t,x)\in\Rpop\times\Rd;
\\ &\tilde U(0,x)=\tilde{u}_{0}(x),\  x\in\Rd;
\end{split}
\eeq	
where $\Rpop=(0,\infty)$.  The wellposedness of \eqref{eq:ksst} in $d\ge2$ is a long-standing open problem.   Consequently, the current KS multidimensional mathematical theory is incomplete (the associated semigroup is not defined everywhere in this case \cite[p.~141]{T}).    

We establish global uniqueness and local existence of time-continuous $\Ltp$-valued, $p\ge1$, solutions for the intimately-linked Burgers incarnation of \eqref{eq:ksst}---the Kuramoto-Sivashinsky-Burgers (KSB) PDE\footnote{Of course, the constants $1/8$ and $4$ in \eqref{eq:ksBur} may easily be changed, and order parameters may be introduced in the differential operator, using scaling parameters in our L-KS kernel as we explicitly detail in \cite{Alksspde}.}:
\beq\lbl{eq:ksBur}
\bsp 
&\displaystyle{\partial_{t} U}=-\tfrac18\lpa\lap^{2}+4\lap\rpa U -\tfrac12\sum_{k=1}^{d}{\pa_{k} U^{2}},\ (t,x)\in\Rpop\times\Rd;
\\ &U(0,x)=\unx,\  x\in\Rd;
\end{split}
\eeq	
in spatial dimensions $d<6p$\footnote{In particular, the unique $\Lt$ solutions exist in spatial dimensions $d=1,2,\ldots,5$.}.  Of interest in its own right due to the KS-Burgers interaction in \eqref{eq:ksBur} and the ubiquitous role the KS and Burgers ingredients play in modeling turbulent behavior, the Burgers formulation of the KS equation \eqref{eq:ksBur} inherits the same KS differential operator from---and is closely related to---the eikonal KS PDE in \eqref{eq:ksst}.    

To prove our result, we use a mixture of ideas and methods from our Brownian-time program \cite{Alksspde}--\cite{Abtp2} and \cite{AX17,AZ01}; analysis (including complex/harmonic); and the stochastic analytic stopping-time argument, adapted to our deterministic setting.   The Brownian-time component of the argument is encapsulated by our L-KS kernel \cite{Alksspde,Abtbs,Abtpspde,Aks,AX17} and the associated novel formulation of \eqref{eq:ksBur}, which play a crucial role throughout the proof.  On one hand, the L-KS kernel was at the heart of the second major direction of our Brownian-time program, leading to the mild formulation and analysis of a large class of other simpler fourth order KS-type PDEs/Stochastic PDEs in multiple spatial dimensions\footnote{The first was the high order/time fractional PDEs connection to the closely-linked Brownian-time stochastic processes that we started in \cite{AZ01,Abtp2}}.  On the other hand, the form of the L-KS kernel provides a bridge between the second order Schr\"odinger equation, via its propagator, and the---quite different setting---of fourth order KS-type equations (see \secref{sec:LKSform} below).    

As in \cite{T}, when $d=1$ the two formulations \eqref{eq:ksst} and \eqref{eq:ksBur} are interchangeable.  This is easily seen by differentiating \eqref{eq:ksst} with respect to  $x$ and letting $U=\pa\tilde U/\pa{x}$ yielding \eqref{eq:ksBur}.  Similarly, when $d>1$, the KS Burgers PDE \eqref{eq:ksBur} is embedded in---and is a major driver of the regularity of---\eqref{eq:ksst}. To see this, simply apply the directional derivative in the direction of the $d$-dimensional vector $v=\frac1d\lang1,1,\ldots,1\rang$ to \eqref{eq:ksst} to obtain
\beq\lbl{eq:dirder}
\bsp
0&=\lpa{\partial_{t}}+ \tfrac18\lpa\lap^{2}+4\lap\rpa\rpa \lpa\tfrac1d\sum_{\l=1}^{d}\partial_{\l}\tilde U\rpa+\tfrac1{2d}\sum_{\l=1}^{d}\partial_{\l}\sum_{k=1}^{d}\lpa{\partial_{k}\tilde U} \rpa^{2}
\\&=\lpa{\partial_{t}}+ \tfrac18\lpa\lap^{2}+4\lap\rpa\rpa \lpa\tfrac1d\sum_{\l=1}^{d}{\partial_{\l}\tilde U}\rpa
\\&\quad+\tfrac1{2d}\sum_{\l=1}^{d}{\partial_{\l}}\lbk\lpa\sum_{k=1}^{d}{\partial_{k}\tilde U} \rpa^{2}-
\sum_{\substack{k_{1}\neq k_{2}\\ k_{1},k_{2}\in\{1,\dots,d\}}}\lpa{\partial_{k_{1}}\tilde U }\rpa\lpa{\partial_{k_{2}} \tilde U}\rpa\rbk.
\end{split}
\eeq
Clearly, the first four terms to the right of the second equality in \eqref{eq:dirder} are those of a KS Burgers PDE\footnote{\lbl{fn:fourthtrm} Here the $1/2$ in front of the sum in \eqref{eq:ksBur} is replaced with $d/2$ (a change in constant only in \eqref{eq:ksBur}).} in $U=\sum_{\l=1}^{d}{\partial_{\l}\tilde U}/d$.  Furthermore, the extra fifth---or last---term in \eqref{eq:dirder} is the off-diagonal part of the fourth; and, as such, its  summands are of the same form and order as those of the fourth (each is a product of two first order spatial derivatives of $\tilde U$)\footnote{Note that if $\tilde U$ is a solution to \eqref{eq:ksst} with a permutation-invariant gradient $\uppi(\nabla\tilde U)=\nabla\tilde U$, for every permutation $\uppi$---equivalently ${\pa\tilde U}(t,x)/{\pa x_{i}}={\pa\tilde U}(t,x)/{\pa x_{j}}\ \forall i,j\in\{1,2,\ldots,d\}$---then the directional derivative of \eqref{eq:ksst}, \eqref{eq:dirder}, reduces to \eqref{eq:ksBur} with $U=(\sum_{\l=1}^{d}\frac{\partial\tilde U}{\pa x_{\l}})/d$. }.  
 So, in addition to its own interest, the study of the KS Burgers PDE \eqref{eq:ksBur} is illuminating in gaining valuable insight into \eqref{eq:ksst} in spatial dimensions $d\ge1$.  We will henceforth in this article exclusively focus our analysis and discussion on the KS Burgers PDE \eqref{eq:ksBur}.

 \subsection{The imaginary-Brownian-time-Brownian-angle kernel formulation}\lbl{sec:LKSform} A key fundamental ingredient, and novelty, in our approach to analyze KS equations is our imaginary-Brownian-time-Brownian-angle kernel or L-KS kernel.  We use the L-KS kernel---introduced in \cite{Aks} and utilized in a fundamental way in \cite{Alksspde,Abtpspde,AX17} in the formulation of different KS-type equations---to give a new formulation and to analyze the KS Burgers PDE \eqref{eq:ksBur}.  This L-KS kernel, which is the Brownian average of an angled Schr\"odinger propagator (see \eqref{eq:vepvthLKS} below), allows us to obtain sharp enough estimates that lead to, among other ramifications, the surprising $6p$-th dimensional limit phenomenon of KS Burgers equations.  As we explained in our earlier work \cite{Alksspde,Abtpspde,Aks,AX17}, the L-KS kernel  is the fundamental solution to the L-KS PDE
\beq\lbl{eq:lkspde}
\frac{\partial u}{\partial t}=-{\tfrac18}\lpa\lap+2\rpa^{2}u, \mbox{ with } u(0,x)=\delta(x),
\eeq
where $\delta(x)$ is the usual Dirac delta function, and is given by 
\beq\lbl{eq:vepvthLKS}
\bsp
\KKStxy&=\int_{-\infty}^{0}\df{\e^{\i s} \e^{-|x-y|^2/2\i s}}{{\lpa2\pi \i s \rpa}^{d/2}}\ptsz ds+\int_{0}^{\infty}\df{\e^{\i s} \e^{-|x-y|^2/2\i s}}{{\lpa2\pi \i s \rpa}^{d/2}}\ptsz ds
\\&=(2\pi)^{-d}\int_{\Rd}\e^{-\frac{t}{8}\left( -2+\lab\xi\rab^{2} \right) ^{2}}\e^{\i \langle\xi, x-y\rangle}d\xi
\\&= (2\pi)^{-d}\int_{\Rd}\e^{-\frac{t}{8}\left( -2+\lab\xi\rab^{2} \right) ^{2}}\cos\lpa{\langle\xi, x-y\rangle}\rpa d\xi,
\end{split}
\eeq
where $\i=\sqrt{-1}$ and where $\ptsz=\e^{-({s^{2}}/{2 t})}/\sqrt{2\pi t}$ is the density of a one-dimensional Brownian motion starting at $0$.  To obtain the last two equalities in \eqref{eq:vepvthLKS}, we inverted the spatial Fourier transform of $\KKStx=\KKStzx$  (see \cite[Lemma 2.1]{Alksspde} and \lemref{lm:basicfacts}, equation \eqref{eq:LKSF}, below).  Since 
${ \e^{-|x-y|^2/2\i s}}/{{\lpa2\pi \i s \rpa}^{d/2}}$
is the famous Schr\"odinger propagator, we call the first representation---given by the first equality in \eqref{eq:vepvthLKS}---the Brownian average of the angled Schr\"odinger propagator form of the L-KS kernel.  We term the second expression, given by the last two equalities in \eqref{eq:vepvthLKS}, the Fourier form of our L-KS kernel.  The L-KS kernel may also be thought of probabilistically as the ``density\footnote{\lbl{fn:density} Of course, as we detailed in \cite{Alksspde,Aks}, the L-KS kernel is \emph{not} a standard probability density function.}'' of our imaginary-Brownian-time-Brownian-angle process (\cite{Alksspde,Abtpspde,Aks})
\begin{equation}
\X_{B}^{x}(t)\eqdef\begin{cases}  X^{x}(\i B(t)) \exp\left(\i B(t)\right),  & B(t)\ge0;\cr
\i X^{-\i x}(-\i B(t))\exp\left(\i B(t)\right),  & B(t)<0;
\end{cases}
\label{KSprocess}
\end{equation}
where the process $X^{x}$ is an $\Rd$-valued Brownian motion (BM) starting from $x\in\Rd$, $\Xix$ is an independent $i\Rd$-valued BM starting at $-ix$ (so that $i\Xix$ starts at $x$), and both are independent of the inner standard $\R$-valued Brownian motion $B$ starting from $0$.    The clock of the outer Brownian motions $\Xx$ and $\Xix$ is replaced by a positive imaginary Brownian time; and the angle of $\X_{B}^{x}$ in the complex plane is the Brownian motion $B$ \cite{Alksspde,Aks}.  We think of the imaginary-time processes $\{\Xx(is),s\ge0\}$ and $\{i\Xix(-is),s\le0\}$ as having the same complex Gaussian distribution on $\Rd$ with the corresponding complex distributional density (or Schr\"odinger propagator)
$${K}^{\text{\tiny{\sc{SP}}}^d}_{\i s;x,y}=\frac{1}{({2\pi i s})^{d/2}}e^{-|x-y|^2/2is}.$$
The L-KS kernel in \eqref{eq:vepvthLKS} is closely tied to the densities of  Brownian-time processes, like the Brownian-time Brownian motion, and its iterates \cite{Atfhosie,Abtbmsie,Abtp2,AZ01,AX17}---which are the fundamental solutions of time-fractional PDEs and higher order PDEs with memory---as detailed in  \cite{Alksspde,Abtbs,Abtpspde,Aks,AX17} (see also \lemref{lm:basicfacts} and \lemref{lm:btbmlksfest} below and the discussions right before).  This underlying Brownian-time construction simultaneously (1) lends a probabilistic flavor to our treatment of the KS equation here and (2) provides a template for handling time-fractional Burgers PDEs as in our upcoming works (\cite{AT} and beyond).  Moreover, our L-KS formulation---valid for all dimensions simultaneously---makes the KS Burgers equation amenable to a wide spectrum of dimension-dependent investigations, in both the deterministic and stochastic PDE settings.  We explore various other aspects of KS Burgers solutions and their behavior in planned future articles.

To prepare our PDE \eqref{eq:ksBur} for the L-KS kernel formulation, we rewrite it slightly as
\beq\lbl{eq:lkspdea}
\bsp
&{\partial_{t} U}=-{\tfrac18}\lpa\lap+2\rpa^{2}U+{\tf1{2}}U -{\tfrac12}\sum_{k=1}^{d}{\pa_{k} U^{2}}
\\&U(0,x)=\unx.
\end{split}
\eeq	
We then formulate \eqref{eq:lkspdea} (equivalently \eqref{eq:ksBur}) in the mild L-KS formulation:
\beq\lbl{eq:ibtbapsol}
\bsp
U(t,x)&=\intrd\KKStxy \uny dy+\tfrac1{2}\intrdzt \KKStsxy U(s,y)dsdy
\\&+\tfrac12\sum_{\l=1}^{d}\int_{\Rd}\int_{0}^{t}\frac{\pa \KKStsxy}{\pa y_{\l}}U^{2}(s,y)dsdy.
\end{split}
\eeq

\brm\lbl{rem:basicass} 
In our results, and unless explicitly stated, all solutions to the KSB PDE \eqref{eq:ksBur} are in the L-KS kernel mild sense \eqref{eq:ibtbapsol}.   Also, we assume throughout that $T>0$ is a fixed time that is otherwise arbitrary; and that $C$ is a generic constant that may be different from one line to the next.      
\erm
\bnt[Function spaces and related notations]
In the statement of our result we let 
$$\Lp(\Rd;\R):=\lbr f:\Rd\to\R;|f|_{p}^{p}=\int_{\Rd}\lab f(x)\rab^{p}dx<\infty\rbr,$$ 
be the usual $\Lp$ space, and we denote by $\C\lpa[0,T];\Lp\lpa\Rd;\R\rpa\rpa$ the set of continuous functions $\{u(t);0\le t\le T\}$ that are $\Lp\lpa\Rd;\R\rpa$-valued. We will also alternate between the notations $\lab u(t)\rab_{p}$ and $\lab u(t,\cdot)\rab_{p}$ whenever beneficial for the presentation flow.  For a more comprehensive list of acronyms and notation used in this paper the reader is referred to \appref{sec:acrnot}.
\ent
\subsection{The main result and brief proof sketch}\lbl{sec:mainsketch}                                                                                                                                                                                                                                                                                                                       We now turn to our main existence, uniqueness, and regularity result of this article for \eqref{eq:ksBur}.    We call any $U$ satisfying \eqref{eq:ibtbapsol} an L-KS solution to \eqref{eq:ksBur}.
\bfr
\bthm[Local existence and global uniqueness of solutions to the KS Burgers PDE \eqref{eq:ksBur} for $d<6p$]\lbl{thm:exunreg}
Assume that $T>0$ and $p\ge1$ are arbitrary and fixed.  If $\un\in \Ltp\lpa\Rd;\R\rpa$ then there exist a time $\tau\in(0,T]$ and a unique L-KS solution $U\in\C\lpa[0,\tau];\Ltp\lpa\Rd;\R\rpa\rpa$ to the KSB PDE \eqref{eq:ksBur} for $d<6p$.  Moreover,  if existence holds for the KSB PDE for all $t\in[0,T]$ then so does uniqueness.
\ethm
\efr
\brm
\ben
\item  Our result says that time-continuous $\Lt$-valued local solutions (the special case $p=1$) exist and are unique in spatial dimensions $d=1,2,3,4,5$. 
\item With minor adaptations of our arguments, we can readily obtain the added weak $\Ltp$-stability in \appref{app:spatialcont} (see \thmref{thm:wstab})
 and the spatial continuity for our solution when $\un\in \Ltp\lpa\Rd;\R\rpa\cap\C\lpa\Rd;\R\rpa$  (see \appref{app:spatialcont} for additional brief remarks).  We carry out a detailed analysis of the exact modulus of continuity of our solution, including the stochastic PDE case in a future article (see \cite{AX17} for the flavor of these types of results in the simpler linear SPDEs case).  
 \item The global wellposedness of the KSB equation in this article is treated in our upcoming article \cite{Aglobks}.
\een
\lbl{rem:spatialcont}
\erm
Naturally, our kernel formulation \eqref{eq:ibtbapsol}, without the $U/2$ term, and argument also simultaneously establish---as a bonus---the same results in \thmref{thm:exunreg} for the slight variant of \eqref{eq:ksBur}:
\beq\lbl{eq:kspdeavar}
\bsp
&{\partial_{t} U}=-{\tfrac18}\lpa\lap+2\rpa^{2}U -{\tfrac12}\sum_{k=1}^{d}{\pa_{k} U^{2}}
\\&U(0,x)=\unx.
\end{split}
\eeq	
We state this fact as the following theorem.
\bfr
\bthm[Existence and uniqueness for the KSB-variant in $d<6p$]\lbl{thm:exunreg2}
Under the assumptions of \thmref{thm:exunreg}, the same conclusions of \thmref{thm:exunreg} hold for the KSB variant \eqref{eq:kspdeavar}.
\ethm\efr
We now briefly give the structure of the rest of the paper together with the main highlights of the proof.  The proof proceeds in two main stages, which we now summarize.

We start in \secref{sec:convest} by obtaining fundamental estimates for the moduli of $\KKStx$ and its temporal and spatial derivatives---along with previously obtained Fourier transform for $\KKStx$ (in \cite[Lemma 2.1]{Alksspde})---we use them to derive sufficiently sharp, dimension-dependent, $\Ltp$ convolution estimates for the operators
\beq\lbl{eq:convlks}
\bsp
({\rm i})\ (\sKKS_{0} z)(t,x)&= \int_{\R^d}\KKStxy z(y)dy, \mbox{ with }(\sKKS_{0} z)(0,x)=z(x),\\
({\rm ii})\  (\sKKS u)(t,x)&= \int_0^t\int_{\R^d}\KKStrxy u(r,y)dy dr, \mbox{ and }\\
({\rm iii})\  (\sKKSl v)(t,x)&=\int_0^t\int_{\R^d}\frac{\partial \KKStrxy}{\partial y_\l}v(r,y)dy dr,\ \l=1,\ldots, d,
\end{split}
\eeq
for $(t,x)\in[0,T]\times\Rd$ and functions $z\in\Ltp\lpa\Rd;\R\rpa$, $u\in\Lgamo\lpa[0,T];\Lp\lpa\Rd;\R\rpa\rpa$ and $v\in\Lgamt\lpa[0,T];\Lp\lpa\Rd;\R\rpa\rpa$, for $p\ge1$ and suitable $\gamo,\gamt>1$.  Beside playing a crucial role in our analysis of \eqref{eq:ksBur}, these estimates reveal the $6p$ dimensionality limit of the KS Burgers PDE \eqref{eq:ksBur}. 

Next, in \secref{sec:local}, we construct the $\Ltp$-valued local solution to our KS Burgers equation \eqref{eq:ksBur} via the $N$-ball approximating sequence of the L-KS formulation \eqref{eq:ibtbapsol}.  For each $N\in\N$, we map each $U$ on the right hand side of \eqref{eq:ibtbapsol} into the $\Ltp$ ball centered at the origin with radius $N$, $\B_{N}(0)\subset\Ltwo(\Rd;\R)$ (see \eqref{eq:sinks} below).  Then, we use the convolution estimates from \secref{sec:convest} together with a fixed point lemma to show the existence of a unique solution $U_{N}\in\C\lpa[0,T];\Ltp\lpa\Rd;\R\rpa\rpa$, to the $N$-th approximation, satisfying $\sup_{t\in[0,T]}|U_{N}(t)|_{2p}<\infty$ for each $N\in\N$.  Using a stopping-time-flavored argument \`a la stochastic analysis (e.g., \cite[equations (3.15) and (2.3), respectively]{Acom,Acom1} and also \cite{DapZ}), we extract a local L-KS solution to \eqref{eq:ibtbapsol}, $U\in\C\lpa[0,\tau_{\infty});\Ltp\lpa\Rd;\R\rpa\rpa$, for some time $0<\tau_{\infty}\le T$.  This local solution $U$ to \eqref{eq:ksBur} may be thought of as obtained by glueing the path of $U_{N}$ onto that of $U_{N-1}$ for all $N\ge2$.  That argument also establishes uniqueness for \eqref{eq:ksBur}.

\section{Fundamental and convolution estimates for the L-KS kernel  and its derivatives}\lbl{sec:convest}

We first establish some $\L^{2p}$ estimates for the convolution operators with the L-KS kernel and with its spatial derivatives in \eqref{eq:convlks}.  We take up the rougher operator $\sKKSl$ first since it is at the intersection of the KS differential operator---through the kernel---and the Burgers nonlinearity in our mild L-KS formulation \eqref{eq:ibtbapsol}.  We then prove corresponding results for the simpler $\sKKS$ and $\sKKS_{0}$.  

\subsection{Fundamental Brownian-time kernel estimates and Fourier transform}  To decode sufficient information from the L-KS kernel, we need a fundamental lemma followed by convolution lemmas used in the analysis of $\sKKS$ and $\sKKSl$.  The estimates in \lemref{lm:basicfacts}, \eqref{eq:lksfundbdskertimeder}, below were intuited from the intimate relation between the L-KS process \cite{Abtpspde,Aks} and its Brownian-time Brownian motion sibling \cite{Abtp2,AZ01} and their kernels (see also \cite{Alksspde,AX17} for more connections and related properties).  This cozy relation between the two was a major motivation for---and was pointed to repeatedly in---our work \cite{AX17}--\cite{AZ01}. A brief discussion is included in \appref{app:briefdisc} (see in particular the corresponding estimates for the BTBM density, given in \eqref{eq:btbmk}, in \lemref{lm:btbmlksfest} equation \eqref{eq:btbmbds}, which are identical modulo constants).
\blm[Fundamental L-KS lemma: bounds and Fourier transform of $\KKStx$]\lbl{lm:basicfacts}
The spatial Fourier transform of  $\KKStx$ is given by\footnote{\lbl{fn;FTdef}We are using the unitary Fourier transform definition
$$\hat{f}(\xi)=\lpa2\pi\rpa^{-\frac d2}\int_{\Rd}f(x)\e^{-\i\lang x,\xi\rang}dx$$
and the notation $\pa^{2}_{t,\l}$ for $\pa^{2}/\pa t\pa y_{\l}$.}
\beq\lbl{eq:LKSF}
\FKKStxi=\lpa2\pi\rpa^{-\frac d2}\e^{-\frac{t}{8}\left(-2+\lab\xi\rab^{2} \right) ^{2}};\ t>0,\xi\in\Rd, d\in\N.
\eeq
Fix an arbitrary $T>0$, let $t\in(0,T]$, and $\l\in\{1,\cdots,d\}$.  The following estimates hold for the L-KS kernel $\KKStxy:\footnote{The bound in \eqref{eq:lksfundbdskertimeder} (a) is simply a rescaled Brownian-time Brownian motion density (see \eqref{eq:btbmk}).}$
\beq\lbl{eq:lksfundbdskertimeder}
\bsp
&(a) \lab\KKStxy\rab\le C_{1}\int_{0}^{\infty}\frac{\e^{-{c_{1}|x-y|^2}/{s}}}{s^{{d}/{2}}} \frac{\e^{-\tfrac{c_{2}s^2}{t}}}{{t}^{1/2}} d{s}, \\
&(b) \lab{\partial_{t}\KKStxy}\rab\le C_{2} \int_{0}^{\infty}\frac{\e^{-{c_{3}|x-y|^2}/{s}}}{s^{{d}/{2}}}\frac{\e^{-\frac{c_{4}s^2}{t}}}{t^{3/2}} ds,\\
&(c) \lab{\pa_{\l}\KKStxy}\rab\le {C_{3}}\int_{0}^{\infty}\frac{e^{-{{c}_{5}|x-y|^2}/{s}}}{s^{{
(d+1)}/{2}}} \frac{e^{-\frac{{c}_{6} s^2}{t}}}{t^{1/2}} d{s},
\\&(d) \lab{\partial_{t,\l}^{2}\KKStxy}\rab\le C_{4}\int_{0}^{\infty}\frac{e^{-\frac{c_{7}|x-y|^2}{s}}}{s^{({d+1})/{2}}}\frac{e^{-\frac{c_{8}s^2}{t}}}{t^{3/2}} ds, 
\end{split}
\eeq
and hence
\beq\lbl{eq:lksL1bdskertimeder}
\bsp
&(a) \lab\KKStdot\rab_{1}\le\tilde C_{1},
\\&(b)\lab{\partial_{t}\KKStdot}\rab_{1}\le\tilde{C}_{2} t^{-1},
\\&(c) \left|\partial_{\l} \KKStdot\right|_1\le\tilde{C}_{3}t^{-1/4},
\\&(d) \left|\partial^{2}_{t,\l} \KKStdot\right|_1\le\tilde{C}_{4}t^{-5/4},
\end{split}
\eeq
for  $(t,x,y)\in(0,T]\times\Rd\times\Rd$ and for some constants $C_{i},\tilde C_{i},c_{j}\in(0,\infty)$, $i=1,2,3$ and $j=1,2,3,4,5,6,7,8$, that may only depend on $T$ and $d$.  
\elm
\bpfs{Proof of \lemref{lm:basicfacts}}  We refer the reader to \cite[Lemma 2.1]{Alksspde} for the proof of \eqref{eq:LKSF}. 

Moving to the proof of \eqref{eq:lksfundbdskertimeder}, we recall that
\beqs
\bsp
\KKStxy&=\int_{\R\setminus\{0\}}\df{\e^{\i s} \e^{-|x-y|^2/2\i s}}{{\lpa2\pi \i s \rpa}^{d/2}}\frac{\e^{-s^{2}/2t}}{\sqrt{2\pi t}} ds.
\end{split}
\eeqs
We will show the estimates for the integral over $(0,\infty)$:
\beq\lbl{eq:fundestpos}
\sJ_{+}(t,|x-y|)=\int_{0}^{\infty}\df{\e^{\i s} \e^{-|x-y|^2/2\i s}}{{\lpa2\pi \i s \rpa}^{d/2}}\frac{\e^{-s^{2}/2t}}{\sqrt{2\pi t}} ds.
\eeq
The argument for the integral over $(-\infty,0)$, $\sJ_{-}(t,|x-y|)$, is very similar and will not be repeated.  


Let ${\mathbb L}_{\epsilon}$ be the line obtained by rotating the positive real line $\Rp$ clockwise
\bfg[h]
\includegraphics{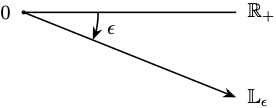}
\ig
\begin{tikzpicture}[scale=.75]
\tkzDefPoint(0,0){O}
\tkzDefPoint(5,0){R}
\tkzDefPoint(5,-2){L}
 
 \draw[thick] (0,0) -- (5,0)
 node[right=2pt]{$\Rp$};
 
 \fill[] (0,0) circle [radius=0.05]
  node[left=3pt]{$0$};
  
  \draw[-Stealth, thick] (0,0) -- (5,-2)
node[right=2pt]{ ${\mathbb L}_{\epsilon}$};

\tkzMarkAngle[size=1.75, arrows=Stealth-, thick](L,O,R)
\tkzLabelAngle[pos=2.1](L,O,R){$\epsilon$}
\end{tikzpicture}
\fi
\caption{}
\label{fig:Figo}
\efg
 in the complex plane $\Cmpl$ by an angle $\epsilon$ (\figref{fig:Figo}). 
Then, $\sJ_{+}$ may be regularized as 
\beqs
\sJ_{+}(t,|x-y|)=\lim_{\epsilon\rightarrow0}\int_{{\mathbb L}_{\epsilon}}\df{\e^{\i s} \e^{-|x-y|^2/2\i s}}{{\lpa2\pi \i s \rpa}^{d/2}}\frac{\e^{-s^{2}/2t}}{\sqrt{2\pi t}} ds.
\eeqs

\begin{figure}[h]
   \centering
\includegraphics{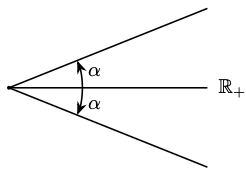}
\ig
\begin{tikzpicture}[scale=0.7]

\tkzDefPoint(0,0){O}
\tkzDefPoint(5,0){R}
\tkzDefPoint(5,-2){D}
\tkzDefPoint(5,2){U}

\draw[thick] (0,0) -- (5,0) 
node[right=2pt]{$\Rp$};

\draw[thick] (0,0) -- (5,2);
\draw[thick] (0,0) -- (5,-2);
\tkzMarkAngle[size=1.85, arrows=-Stealth, thick](R,O,U)
\tkzLabelAngle[pos=2.2](U,O,R){$\alpha$}

\tkzMarkAngle[size=1.85, arrows=Stealth-, thick](D,O,R)
\tkzLabelAngle[pos=2.2](D,O,R){$\alpha$}

\fill[] (0,0) circle [radius=0.05]
  node[left=3pt]{};
\end{tikzpicture}
\fi
\caption{}
\label{fig:Figt}
\end{figure}

Moreover, the integral over ${\mathbb L}_{\epsilon}$ does not depend on $\epsilon>0$ and, for small $\epsilon>0$,  is well defined for complex $t$ inside the cone whose vertex is at the complex plane origin with direction along $\Rp$---and to the right---and aperture $\alpha$ (the region in \figref{fig:Figt} containing $\Rp$ and bounded above and below  by the two lines starting at the complex plane origin and making a fixed angle $\alpha$ above and below $\Rp$), 
provided $\alpha+2\epsilon<\pi/2$.
Take $\alpha=\epsilon=\pi/8$.  Then, on ${\mathbb L}_{\epsilon}$, we have
\beq\lbl{eq:comptoreal}
\bsp
\Re(\i s)\le|s|,\ \Re\lpa1/{\i s}\rpa\ge\sin(\pi/8)\cdot\frac1{|s|}, \mbox{ and }\Re(s^{2}/2t)\ge\sin(\pi/8)\cdot\frac{|s|^{2}}{2|t|}.
\end{split}
\eeq
Thus,
\beq\lbl{eq:J+bd0}
\lab\sJ_{+}(t,|x-y|)\rab\le\frac{|t|^{-1/2}}{\sqrt{2\pi}}\int_{0}^{\infty}\frac{\e^{\sigma}\e^{-\frac{|x-y|^2}{2\sigma}\sin(\pi/8)}\e^{-\frac{\sigma^2}{2|t|}\sin(\pi/8)}}{(2\pi\sigma)^{{d}/{2}}} d\sigma.
\eeq
We now note that 
\beq\lbl{eq:expineq}
\bc
\e^{\sigma-{\frac{\sigma^2\sin(\pi/8)}{2|t|}}}\le \e^{-{\frac{\sigma^2\sin(\pi/8)}{4|t|}}},&\mbox{ for }\sigma\ge4|t|/\sin(\pi/8);\\
\e^{\sigma}\le \e^{\frac{4|t|}{\sin(\pi/8)}},&\mbox{ for }\sigma<4|t|/\sin(\pi/8).
\ec
\eeq
We then easily obtain the bound
\beq\lbl{eq:J+bdcomp}
\bsp
\lab\sJ_{+}(t,|x-y|)\rab\le\frac{\e^{\frac4{\sin(\pi/8)}|t|}|t|^{-1/2}}{\sqrt{2\pi}}\int_{0}^{\infty}\frac{e^{-\frac{|x-y|^2}{2\sigma}\frac{\sin(\pi/8)}{4}}e^{-\frac{\sigma^2}{2|t|}\frac{\sin(\pi/8)}{4}}}{(2\pi\sigma)^{{d}/{2}}} d\sigma.
\end{split}
\eeq
Finally, for any fixed $t\in(0,T]$ use the bound in \eqref{eq:J+bdcomp} in the complex disk centered at $t$ with radius $t\sin(\pi/8)$ (\figref{fig:Figth}) 
\bfg[h]
\centering
\includegraphics{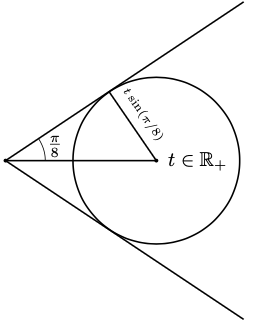}
\ig
\begin{tikzpicture}[scale=0.7]
\coordinate (P) at (0,0);
\coordinate (Q) at (6,4);
\coordinate (R) at (6,-4);
\coordinate (C) at (3.8,0);
\coordinate(D) at (4,0);
\coordinate(F) at (5.9,2); 
\coordinate(U) at (2.6,1.75);

\draw[thick] (P)--(Q);  
\draw[thick] (P)--(R);

\draw[thick] (P)--(C); 
\draw[thick] (C)--(U);
\fill[] (C) circle [radius=0.05];
\fill[] (P) circle [radius=0.05];

\tkzLabelSegment[rotate=305,above=2pt](C,U){\tiny{$t\sin({\pi}/8)$}} 
\tkzLabelSegment[right=0.5pt](C,D){{$t\in \mathbb{R}_{+}$} }
\tkzMarkAngle[mark=none, size=1.cm](C,P,Q)
\tkzLabelAngle[pos=1.3](C,P,Q){$\frac{\pi}8$}

\draw[thick] (C) circle [radius=2.1];
\end{tikzpicture}
\fi
\caption{}
\label{fig:Figth}
\efg

together with Cauchy's estimate to get
\beq\lbl{eq:J+bdreal}
\bsp
&(a) \lab\sJ_{+}(t,|x-y|)\rab\le C_{T,d}t^{-1/2}\int_{0}^{\infty}\frac{e^{-\frac{|x-y|^2}{2\sigma}\frac{\sin(\pi/8)}{4}}e^{-\frac{\sigma^2}{2t}\frac{\sin(\pi/8)}{4}}}{\sigma^{{d}/{2}}} d\sigma,
\\&(b)\lab{\pa_{t}}\sJ_{+}(t,|x-y|)\rab\le C_{T,d}t^{-3/2}\int_{0}^{\infty}\frac{e^{-\frac{|x-y|^2}{2\sigma}\frac{\sin(\pi/8)}{4}}e^{-\frac{\sigma^2}{2t}\frac{\sin(\pi/8)}{4}}}{\sigma^{{d}/{2}}} d\sigma,
\end{split}
\eeq
where the constant $C_{T,d}\in(0,\infty)$ depends on $T$ and $d$.  
Similarly we obtain the desired bounds for $\lab\sJ_{-}(t,|x-y|)\rab$ and $\lab{\pa_{t}}\sJ_{-}(t,|x-y|)\rab$.  This establishes \eqref{eq:lksfundbdskertimeder} parts (a) and (b).

Replacing $\KKStxy$ with
\beqs
\pa_{\l}\KKStxy=\int_{\R\setminus\{0\}}\df{\e^{\i s} (x_{\l}-y_{\l})\e^{-|x-y|^2/2\i s}}{{\lpa2\pi \rpa}^{d/2} (\i s)^{\tf d2+1}}\frac{\e^{-s^{2}/2t}}{\sqrt{2\pi t}} ds,
\eeqs
and applying the same argument above with only minor changes, we obtain
\beq\lbl{eq:derker}
\bsp
\lab\pa_{\l}\KKStxy\rab\le C_{T,d}\int_{0}^{\infty}\lab(x_{\l}-y_{\l})\frac{e^{-c{|x-y|^2}}}{s^{\tf d2+1}}\rab\frac{\e^{-\tilde cs^{2}/t}}{\sqrt{t}}ds,
\end{split}
\eeq
for some constants $c,\tilde c,C_{T,d}\in(0,\infty)$, that may depend on $T$ and $d$.  The estimate in \eqref{eq:lksfundbdskertimeder} (c) now readily follows upon using \eqref{eq:derker} together with	
the following standard estimate for the derivative of the Brownian motion density (or heat kernel):
\beq\lbl{eq:stndrdspcederest}
\lab(x_{\l}-y_{\l})\frac{e^{-{|x-y|^2}/{2s}}}{(2\pi )^{{
d}/{2}}s^{\tf d2+1}}\rab=\lab\frac{\pa}{\pa y_{\l}}\frac{e^{-{|x-y|^2}/{2s}}}{(2\pi s)^{{
d}/{2}}}\rab\le C\frac{e^{-{c|x-y|^2}/{s}}}{s^{{(d+1)}/{2}}}, 
\eeq
for every  $(s,x,y)\in(0,\infty)\times\Rd\times\Rd, \l\in\{1,\cdots,d\}$, for some constants $c,C\in(0,\infty)$ that may depend on $d$;  and, as above, Cauchy's estimate then gives us  \eqref{eq:lksfundbdskertimeder} (d).  
The proof is complete since \eqref{eq:lksL1bdskertimeder} follows immediately from \eqref{eq:lksfundbdskertimeder} upon integration over $\Rd$.
\epfs
\back
I would like to thank Fedja Nazarov for his suggestions about the  proof of \eqref{eq:lksfundbdskertimeder} in \lemref{lm:basicfacts}, which we adopted above in lieu of our original longer argument.  I'd also like to thank Hala Eweiss Allouba for her help with the figures in the proof and, more importantly, for her invaluable love and support.  
\eack
\subsection{Convolution with the L-KS kernel derivative: the $6p$-th dimensional barrier}
We now establish regularity estimates for the operator $\sKKSl$ in \eqref{eq:convlks} linked to the Burgers terms sum of our KS equation (second on the right hand side  of \eqref{eq:ksBur} or third in \eqref{eq:lkspdea} and its mild formulation \eqref{eq:ibtbapsol}).
\blm[Estimates for the convolution operator $\sKKSl$]  \lbl{lm:Burgersfundest} Fix arbitrary $T>0$ and $p\ge1$, and let $d<6p$.  For all  {$\gamma> {8p}/(6p-d)$} the operator $\sKKSl$ is bounded from $\Lgam\lpa[0,T];\Lp\lpa\Rd;\R\rpa\rpa$ into $\C\lpa[0,T];\Ltp\lpa\Rd;\R\rpa\rpa$, and it satisfies the following estimates: 
\beq\lbl{eq:cnvlem1est1}
\bsp
&(a) \left|(\sKKSl v)(t)\right|_{2p}\leq C\int_0^t (t-r)^{-\frac{2p+d}{8p}}|v(r)|_{p} dr\leq C\left(\int_0^t|v(r)|_p^\gamma dr\right)^{1/\gamma},
\\&(b) {\left|(\sKKSl v)(t)-(\sKKSl v)(s)\right|_{2p}\leq C(t-s)^\varrho \left(\int_0^t|v(r)|_p^\gamma dr\right)^{1/\gamma}},
\end{split}
\eeq
for  some constant $C\in(0,\infty)$---that may depend on $T$, $d$, and $p$---for $\varrho<\frac{6p-d}{8p}-\frac{1}{\gamma}$, and for $0\le s<t\le T$.  
\elm
\brm\lbl{rem:Fivedimlimit}
The bounds on the time exponents $-(2p+d)/8p$, and hence $\gamma>8p/(6p-d)$, in our estimates above is the reason our solution is limited to  spatial dimensions $d<6p$.  
\erm
\bpf  We start with the proof of \eqref{eq:cnvlem1est1} part (a).  Using Minkowski's and Young's inequalities and denoting the spatial convolution by $*$ and the spatial partial derivative in the $\l$-th spatial variable by $\pa_{\l}$,  give us 
\beq\lbl{eq:L2pkerder}
\bsp
\left|(\sKKSl v)(t)\right|_{2p}
&=\lbk\int_{\Rd}\lab\int_0^t\lpa{\partial_{\l}}{\KKStrdot}*v(r,\cdot)\rpa(x) dr\rab^{2p}dx\rbk^{1/2p}
\\&\leq  \int_0^t\left|\lab\partial_{\l} \KKStrdot\rab*\lab v(r,\cdot)\rab\right|_{2p}dr
\\&\leq \int_0^t \left|\partial_{\l} \KKStrdot\right|_q\left|v(r,\cdot)\right|_p dr, 
\end{split}
\eeq
where $q=\frac{2p}{2p-1}$.  To finish the proof of \eqref{eq:cnvlem1est1} (a), we need to estimate $|\partial_{\l} \KKStrdot|_q$ for $p\ge1$ (equivalently $1<q\le2$).  We carry out this in three steps:
\ben
\item we begin with the case $p_{0}=1$ ($q_{0}=2$), 
\item then we obtain the corresponding estimate for $p> p_{1,d}>p_{0}$, for some suitable dimension-dependent $p_{1,d}$, and 
\item then we appeal to the standard log-convexity of $\Lp$-norms to interpolate through the gap between the $p_{0}$ and $p_{1,d}$ values (or their equivalent $q$ values) in the estimates for $|\partial_{\l} \KKStrdot|_q$.
\een

For the case $p_{0}=1$ ($q_{0}=2$), let $\S\subset\Rd$ be given by
\beq\lbl{eq:thesetS}
\S=\lbr\xi\in\Rd;|\xi|\in\lbk0,\sqrt{4-2\sqrt2}\ \rbk\cup\lbk\sqrt{4+2\sqrt2},\infty\rpa\rbr,
\eeq 
and observe that, for $0<t\le T$, we have
\beqs
\bc
\e^{-\frac t4(-2+|\xi|^{2})^{2}}\le\e^{-\frac t8|\xi|^{4}}, &\xi\in\S;\cr
\e^{-\frac t4(-2+|\xi|^{2})^{2}}\le C_{T}\e^{-\frac t8|\xi|^{4}}, &\xi\in\Rd\setminus\S;
\ec
\eeqs
where $C_{T}\in(1,\infty)$ depends only on $T$.
Then, by the Parseval-Plancherel identity and the Fourier transform of the L-KS kernel, in \lemref{lm:basicfacts} equation \eqref{eq:LKSF}, we have that 
\beq\lbl{eq:Parskerderl2}
\bsp
\left|\partial_{\l} \KKStrdot\right|_2&=(2\pi)^{-\frac d2}\lbk\int_{\R^d}\xi_{\l}^{2} \e^{- \frac{(t-r)}{4}(-2 + |\xi|^{2})^{2}} d\xi\rbk^{\frac12}
\\&\le (2\pi)^{-\frac d2}\lbk\int_{\S}\xi_{\l}^{2} \e^{- \frac{(t-r)}{8}|\xi|^{4}} d\xi+C_{T}\int_{\Rd\setminus\S}\xi_{\l}^{2} \e^{- \frac{(t-r)}{8}|\xi|^{4}} d\xi\rbk^{\frac12}
\\&\le (2\pi)^{-\frac d2}\sqrt{C_{T}}\lbk\int_{\R^d}\xi_{\l}^{2} \e^{- \frac{(t-r)}{8}|\xi|^{4}} d\xi\rbk^{\frac12}
\\&=C|t-r|^{-\frac{d+2}{8}},
\end{split}
\eeq
for $0\le r<t\le T$ for all $d\in\N=\{1,2,3, \ldots\}$.  

For the $1<q<2$ ($p>1$) case, we start by using \lemref{lm:basicfacts}, equation \eqref{eq:lksfundbdskertimeder} (c), to derive
\beq\lbl{eq:kerderl2q0}
\bsp
\left|\partial_{\l} \KKStrdot\right|_{q}&\le C\lbk\int_{\Rd}\lab\int_{0}^{\infty}\frac{e^{-{{c}_{5}|x|^2}/{s}}}{s^{{
(d+1)}/{2}}} \ugatmr(ds)\rab^{q}dx\rbk^{1/q},
\end{split}
\eeq
where 
$$\ugatmr(ds)=\frac{e^{-\frac{{c}_{6} s^2}{t-r}}}{\sqrt{{\pi(t-r)}/{4 {c}_{6}}}} d{s}$$ 
is the Gaussian probability measure on $(0,\infty)$.  We will use the elementary integral 
\beq\lbl{eq:elemGauss}
\int_0^\infty s^\alpha \frac{e^{-\frac{{c}_{6} s^2}{t}}}{\sqrt{{\pi t}/{4 {c}_{6}}}} d{s}={\frac {{t}^{-\tf\alpha{2}}{{c_6}}^{\tf\alpha{2}}\Gamma \left( \frac{1-\alpha}{2}
 \right) }{\sqrt {\pi}}}, \mbox{ for }\alpha<1, t>0.
\eeq

We now handle the $d=1$ and the multidimensional cases separately.
 we apply Jensen's inequality to \eqref{eq:kerderl2q0} to obtain
\beq\lbl{eq:kerderl2q1}
\bsp
\left|\partial_{\l} \KKStrdot\right|_{q}&\le C\lbk\int_{0}^{\infty}\int_{\Rd}\lpa\frac{e^{-{{c}_{5}|x|^2}/{s}}}{s^{{
(d+1)}/{2}}}\rpa^{q}dx\, \ugatmr(ds)\rbk^{1/q}
\\&\le C\lbk\int_{0}^{\infty}s^{-(\frac{d+1}{2}q-\frac d2)}\ugatmr(ds)\rbk^{1/q}
\\&= C (t-r)^{-\frac{d(q-1)+q}{4q}} 
= C (t-r)^{-{\frac {2\,p+d}{8p}}},
\end{split}
\eeq
for $1<q<(d+2)/(d+1)$ (equivalently $p>{(d+2)}/{2}$).  In particular, for $d=1$ $1<q<3/2$ (equivalently $p>3/{2}$).

On the other hand, applying Minkowski's inequality to \eqref{eq:kerderl2q0} we get  
\beq\lbl{eq:kerderl2q2}
\bsp
\left|\partial_{\l} \KKStrdot\right|_{q}&\le C\int_{0}^{\infty}\lbk\int_{\Rd}\lpa\frac{e^{-{{c}_{5}|x|^2}/{s}}}{s^{{(d+1)}/{2}}}\rpa^{q}dx\rbk^{1/q}
 \ugatmr(ds)\\&\le C\int_{0}^{\infty}s^{-(\frac{d+1}{2}-\frac {d}{q})}\ugatmr(ds)
\\&= C (t-r)^{-\frac{d(q-1)+q}{4q}} 
= C (t-r)^{-{\frac {2\,p+d}{8p}}},
\end{split}
\eeq
where $q<\frac d{d-1}$ (equivalently $p>\frac d2$) and $d\ge2$.\footnote{Minkowski's inequality yields the desired bound for a wider range of $q$ (and hence $p$) than Jensen's, but since $q<2$ the Minkowski bound $q<\frac d{d-1}$ is useful only for $d\ge2$ and we use the Jensen's inequality bound for $d=1$.  
}

 We have now arrived at the log-convexity interpolation step (3) we outlined above to finish the $|\partial_{\l} \KKStrdot|_{q}$ estimates for the remaining $p$ values  $1=p_0<p\leq p_{1,d}$, where $p_{1,d}=3/2$ for $d=1$ and $p_{1,d}=d/2$ for $d\ge2$ (by the comments after \eqref{eq:kerderl2q1} and \eqref{eq:kerderl2q2}, respectively).  So, for $0<\theta<1$, let $p_{1,d}^{+}>p_{1,d}$ be fixed but arbitrary and let $q_{1,d}^{+}=2p_{1,d}^{+}/(2p_{1,d}^{+}-1)$.  We use the standard relation
\beqs
\bsp
\frac1{p_{\theta}}=\frac{(1-\theta)}{p_{0}}+\frac\theta{p_{1,d}^{+}}\mbox{ and } 
\frac1{q_\theta}=\frac{(1-\theta)}{q_{0}}+\frac\theta{q_{1,d}^{+}} 
\end{split}
\eeqs
and the $\Lp$ norms log-convexity readily gives
\beq\lbl{eq:interpoltn}
\bsp
\left| \partial_{\l} \KKStrdot\right|_{q_\theta} 
&\le\left| \partial_{\l} \KKStrdot\right|_{q_0}^{1-\theta}\left| \partial_{\l} \KKStrdot\right|_{q_{1,d}}^{\theta}
\\&\leq C (t-r)^{-\frac{d(q_\theta-1)+q_\theta}{4q_\theta}}  
\\&= C (t-r)^{-\frac {2p_\theta+d}{8p_\theta} }, 1<p_{\theta}<p_{1,d}^{+}.
\end{split}
\eeq
By equations \eqref{eq:Parskerderl2}, \eqref{eq:kerderl2q1}, \eqref{eq:kerderl2q2}, \eqref{eq:interpoltn}, and the definitions of $p_{1,d}$ and $p_{1,d}^{+}$ above we get that the $\Lq$ bound on 
 $\partial_{\l} \KKStrdot$ is given by 
 \beq\lbl{eq:lqderkergen}
 \lab\KKStrdot\rab_{q}\le C (t-r)^{-\frac{d(q-1)+q}{4q}} 
=C (t-r)^{-{\frac {2\,p+d}{8p}}}, \ \forall p\in[1,\infty),\ d\in\N.  
 \eeq
 Thus, by \eqref{eq:L2pkerder} and H\"older's inequality we obtain  
$$\left|(\sKKSl v)(t)\right|_{2p}\leq C\int_0^t (t-r)^{-{\frac {2p+d}{8p}}}|v(r)|_p dr\leq C\left(\int_0^t|v(r)|_p^\gamma dr\right)^{1/\gamma},$$
for $\gamma >\frac {8p}{6p-d}$, $d<6p$, $p\ge1$, and $t\in(0,T]$.  This establishes \eqref{eq:cnvlem1est1} (a).

We now turn to \eqref{eq:cnvlem1est1} part (b).   For $s<t$, we have
\beq\lbl{eq:T1T2}
\left|(\sKKSl v)(t)-(\sKKSl v)(s)\right|_{2p}\leq \sum_{i=1}^{2} F_{i}^{(p)}(s,t),
\eeq
where
\beqs 
F_{1}^{(p)}(s,t)=\lab \int_s^t \int_{\Rd}  \frac{\partial \KKStrdoty}{\partial y_\l} v(r,y)dy dr\rab_{2p}
\eeqs
and
\beqs 
\bsp
F_{2}^{(p)}(s,t)=\lab \int_0^s\int_{\Rd}  \left( \frac{\partial\KKStrdoty}{\partial y_\l} - \frac{\pa\KKSsrdoty}{\partial y_\l}\right) v(r,y)dy dr\rab_{2p}.
\end{split}
\eeqs
We bound the first term $F_{1}^{(p)}$ using exactly the same argument in \eqref{eq:L2pkerder} leading to \eqref{eq:cnvlem1est1} (a).  So, using Minkowski's and Young's inequalities, our estimate on $|\partial_{\l} \KKStrdot|_{q}$, and then H\"older's inequality we get
\beq\lbl{eq:T1}
\bsp
F_{1}^{(p)}(s,t)&\leq \int_s^t\left|\partial_{\l} \KKStrdot\right|_{q}|v(r)|_p dr
\\&\leq C\int_s^t(t-r)^{-\frac{2p+d}{8p}}|v(r)|_p dr
\\&\leq C\left[\int_s^t (t-r)^{-\frac{2p+d}{8p}\frac{\gamma}{\gamma -1}}dr\right]^{\frac{\gamma-1}{\gamma }}\left[\int_s^t |v(r)|_p^\gamma dr\right]^{1/\gamma}
\\&\leq C|t-s|^{1-\frac{2p+d}{8p}-\frac{1}{\gamma}}\lbk\int_0^t |v(r)|_p^\gamma dr\rbk^{1/\gamma},
\end{split}
\eeq
for $\gamma >\frac {8p}{6p-d}$, $d<6p$, $p\ge1$, and $0\le s<t\le T$.

For  $F_{2}^{(p)}$, we use Minkowski's and Young's inequalities as in \eqref{eq:L2pkerder} to get
\beq\lbl{eq:T2.1inq1}
\bsp
F_{2}^{(p)}(s,t)&\leq\int_0^s \lab\partial_{\l} \KKStrdot-\partial_{\l} \KKSsrdot\right|_q \left|v(r,\cdot)\rab_p dr,
\end{split}
\eeq 
Let $\pa_{\theta}=\pa/\pa{\theta}$.  By Minkowski's inequality, the fundamental theorem of calculus,  \lemref{lm:basicfacts} \eqref{eq:lksfundbdskertimeder} (c) and (d), and  the $|\partial_{\l} \KKStrdot|_{q}$ estimates above yield 
\beq\lbl{eq:lqdiffest}
\bsp
\lab\partial_{\l} \KKStrdot-\partial_{\l} \KKSsrdot\right|_{q}&\le\int_{s}^{t}\lab\pa_{\theta}\partial_{\l} \KKSthrdot\rab_{q}d\theta
\\&\le C\int_{s}^{t}(\theta-r)^{-1}\lab\partial_{\l} \KKSthrdot\rab_{q}d\theta
\\&\le C\int_{s}^{t}(\theta-r)^{-1-\tfrac{2p+d}{8p}}d\theta,
\end{split}
\eeq
for $0<r<s<\theta<t\le T$, $p\ge1$, $d\in\N$, and for a constant $C\in(0,\infty)$ that depends only on $T$. 
 
 Combining \eqref{eq:T2.1inq1} with \eqref{eq:lqdiffest}, and using H\"older's inequality in $\theta$ and then in $r$  we get
\beq\lbl{eq:T2.1inq4}
\bsp
F_{2}^{(p)}(s,t)&\le C\int_0^s \int_{s}^{t}(\theta-r)^{-1-\tfrac{2p+d}{8p}} \left|v(r,\cdot)\rab_p d\theta dr
\\&\le C|t-s|^{\varrho}\int_{0}^{s}(s-r)^{1-\varrho-\tfrac{10p+d}{8p}}\left|v(r,\cdot)\rab_p dr
\\&\le C|t-s|^{\varrho}\lpa\int_{0}^{t}\left|v(r,\cdot)\rab_p^{\gamma} dr\rpa^{1/\gamma},
\end{split}
\eeq 
for $\gamma >\frac {8p}{6p-d}$, $\varrho<1-\frac{2p+d}{8p}-\frac{1}{\gamma}$,  $d<6p$, and $0<s<t\le T$.
Thus, the estimate \eqref{eq:cnvlem1est1} (b) follows from equations \eqref{eq:T1T2}, \eqref{eq:T1}, and  \eqref{eq:T2.1inq4}.
The proof of the lemma is now complete.  
\epf
\subsection{Convolutions with the L-KS kernel}
We now state and prove the estimates for the operators $\sKKS$ and $\sKKS_{0}$ in \eqref{eq:convlks}, which are connected to the second and initial data terms in \eqref{eq:lkspdea} and its mild formulation \eqref{eq:ibtbapsol}).  
\blm[Estimates for the convolution operator $\sKKS$] \lbl{lm:bdsKScnv} Fix an arbitrary $T>0$  and $p\ge1$, and let $d<6p$.  For every $\gamma> 1$ the operator $\sKKS$ is bounded from $\Lgam\lpa[0,T];\Ltp\lpa\Rd;\R\rpa\rpa$ into $\C\lpa[0,T];\Ltp\lpa\Rd;\R\rpa\rpa$, and it satisfies the following estimates: 
\beq\lbl{eq:cnvlem2est1}
\bsp
(a)& \left|(\sKKS u)(t)\right|_{2p} \le C\int\limits_0^t|u(r)|_{2p} dr,
\\ (b)& \left|(\sKKS u)(t)-(\sKKS u)(s)\right|_{2p} \le C(t-s)^\varrho \left(\int\limits_0^t|u(r)|_{2p}^\gamma dr\right)^{1/\gamma},
\end{split}
\eeq
for  some constant $C\in(0,\infty)$---that may depend on $T$, $d$, and $p$---for $\varrho<1-\frac{1}{\gamma}$  and for $0\le s<t\le T$.
  \elm
\bpf  We start with the proof of \eqref{eq:cnvlem2est1} part (a).  Using Minkowski's and Young's inequalities and denoting the spatial convolution by $*$,  give us 
\beq\lbl{eq:L2kerlem2}
\bsp
\left|(\sKKS u)(t)\right|_{2p}
&=\lbk\int_{\Rd}\lab\int_0^t\lpa{\KKStrdot}*u(r,\cdot)\rpa(x) dr\rab^{2p}dx\rbk^{1/2p}
\\&\leq  \int_0^t\left|\lab \KKStrdot\rab*\lab u(r,\cdot)\rab\right|_{2p}dr
\\&\leq \int_0^t \left| \KKStrdot\right|_1 \left|u(r,\cdot)\right|_{2p} dr\le C\int_0^t \left|u(r,\cdot)\right|_{2p} dr, 
\end{split}
\eeq 
where the last inequality follows from \eqref{eq:lksL1bdskertimeder} (a) in \lemref{lm:basicfacts}.  So, \eqref{eq:cnvlem2est1} (a) is established. 

We now handle \eqref{eq:cnvlem2est1} (b) similarly to the corresponding part in \lemref{lm:Burgersfundest}.   For $s<t$, we have
\beq\lbl{eq:T1T2lem2}
\left|\sKKS(u)(t)-\sKKS(u)(s)\right|_{2p}\leq \sum_{i=1}^{2} F^{(p)}_{i}(s,t),
\eeq
where
\beqs 
F^{(p)}_{1}(s,t)=\lab \int_s^t \int_{\Rd}  {\KKStrdoty} u(r,y)dy dr\rab_{2p}
\eeqs
and
\beqs 
\bsp
F^{(p)}_{2}(s,t)=\lab \int_0^s\int_{\Rd}  \left( {\KKStrdoty} - {\KKSsrdoty}\right) u(r,y)dy dr\rab_{2p}.
\end{split}
\eeqs
Bounding $F^{(p)}_{1}$ as in \eqref{eq:L2kerlem2} followed by \lemref{lm:basicfacts} \eqref{eq:lksfundbdskertimeder} (a) and H\"older's inequality
\beq\lbl{eq:F1bdlem2}
\bsp
F^{(p)}_{1}(s,t)&\le \int_s^t\left| \KKStrdot\right|_1\left|u(r,\cdot)\right|_{2p} dr
\\&\le C \int_s^t \left|u(r,\cdot)\right|_{2p} dr\le C(t-s)^{\tf{\gamma-1}{\gamma}}\lpa\int_s^t \left|u(r,\cdot)\right|_{2p}^{\gamma} dr\rpa^{\tf1\gamma}.
\end{split}
\eeq
For $F^{(p)}_{2}$, we proceed via Minkowski's inequality; the fundamental theorem of calculus; the L-KS kernel estimates in \lemref{lm:basicfacts} \eqref{eq:lksfundbdskertimeder} (b) and \eqref{eq:lksL1bdskertimeder} (b);
and Minkowski's, Young's and H\"older's inequalities to get
\beq\lbl{eq:F2bdlem2}
\bsp
&F^{(p)}_{2}(s,t)\le \int_{0}^{s}\lab\int_{s}^{t}\int_{\Rd}\lab\pa_{\theta}\KKSthrdoty\rab\lab u(r,y)\rab dy d\theta\rab_{2p}dr
\\&\le C\int_{0}^{s}\lbk\int_{\Rd}\lab\int_{s}^{t}\int_{\Rd}\int_{0}^{\infty}\frac{e^{-\frac{{c_{1}|x-y|^2}}{\varsigma}}{e^{-\frac{c_{2}\varsigma^2}{\theta-r}}}}{\varsigma^{{d}/{2}}{(\theta-r)^{3/2}}} \lab u(r,y)\rab d\varsigma dy d\theta\rab^{2p}dx\rbk^{1/2p} dr
\\&\le C\int_{0}^{s}\int_{s}^{t}\lab\int_{0}^{\infty}\frac{e^{-\frac{{c_{1}|\cdot|^2}}{\varsigma}}{e^{-\frac{c_{2}\varsigma^2}{\theta-r}}}}{\varsigma^{{d}/{2}}{(\theta-r)^{3/2}}} d\varsigma\rab_{1}\lab u(r,\cdot)\rab_{2p}d\theta dr
\\&\le C\int_{0}^{s}\int_{s}^{t}(\theta-r)^{-1}\lab u(r,\cdot)\rab_{2p}d\theta dr
\\&\le C(t-s)^{\varrho}\int_{0}^{s}(s-r)^{\beta-1}\lab u(r,\cdot)\rab_{2p} dr
\\&\le C(t-s)^{\varrho}T^{1-\frac1\gamma-\varrho}\lbk\int_{0}^{t}\lab u(r,\cdot)\rab_{2p}^{\gamma} dr\rbk^{1/\gamma}
\\&\le C(t-s)^{\varrho}\lbk\int_{0}^{t}\lab u(r,\cdot)\rab_{2p}^{\gamma} dr\rbk^{1/\gamma},
\end{split}
\eeq 
where $\varrho+\beta=1$, provided $\varrho<1-1/\gamma$ and $0<s<t\le T$.  Thus, the estimate \eqref{eq:cnvlem2est1} (b) follows from equations \eqref{eq:T1T2lem2}, \eqref{eq:F1bdlem2}, and  \eqref{eq:F2bdlem2}.

The proof of the lemma is therefore complete.  
\epf 
Using \lemref{lm:basicfacts} and arguing similarly to the proof of \lemref{lm:bdsKScnv}, only simpler, we readily obtain the next result for $\sKKS_{0}$ in \eqref{eq:convlks}. 
\blm[Estimates for the convolution operator $\sKKS_{0}$] \lbl{lm:bdsKSinitcnv} Fix an arbitrary $T>0$, $p\ge1$, and let  $d<6p$.  The operator $\sKKS_{0}$ is bounded from $\Ltp\lpa\Rd;\R\rpa$ into $\C\lpa[0,T];\Ltp\lpa\Rd;\R\rpa\rpa$, and it satisfies the following estimates: 
\beq\lbl{eq:cnvlem3est1}
\bsp
(a)& \left|(\sKKS_{0} z)(t)\right|_{2p}  \le C|z|_{2p}, \mbox{ with } \left|(\sKKS_{0} z)(0)\right|_{2p}=|z|_{2p}, 
\\ (b)& \left|(\sKKS_{0} z)(t)-(\sKKS_{0} z)(s)\right|_{2p}  \le C\lbk\left(\log(t)-\log(s)\right)\nim 1\rbk |z|_{2p} ,
\end{split}
\eeq
for $0<s<t\le T$.
\elm
\bpf
Young's inequality together with \lemref{lm:basicfacts},  \eqref{eq:lksL1bdskertimeder} (a),  give us 
\beq\lbl{eq:unzz}
\left|(\sKKS_{0} z)(t)\right|_{2p} \le\left|\lab\KKStdot\rab*\lab z\rab\right|_{2p} \le\lab\KKStdot\rab_{1}\lab z\rab_{2p}\le C\lab z\rab_{2p} <\infty, 
\eeq
for $0<t\le T$; and the $t=0$ case follows trivially from \eqref{eq:convlks} (i), establishing estimate \eqref{eq:cnvlem3est1} (a).

Arguing as in \eqref{eq:F2bdlem2} above we obtain
\beq\lbl{eq:logest}
\bsp
\left|(\sKKS_{0} z)(t)-(\sKKS_{0} z)(s)\right|_{2p}&\le C\int_{s}^{t}\lab\int_{0}^{\infty}\frac{e^{-\frac{{c_{1}|\cdot|^2}}{\varsigma}}{e^{-\frac{c_{2}\varsigma^2}{\theta}}}}{\varsigma^{{d}/{2}}{\theta^{3/2}}} d\varsigma\rab_{1}d\theta\lab z\rab_{2p}
\\&\le C(\log(t)-\log(s))\lab z\rab_{2p},
\end{split}
\eeq
for $0<s<t\le T$.  The bound in \eqref{eq:cnvlem3est1} (b) is now trivially obtained using \eqref{eq:logest} and \eqref{eq:cnvlem3est1} (a) together with Minkowski's inequality.  We turn to the continuity at time $0$.  In light of \lemref{lm:basicfacts} \eqref{eq:lksfundbdskertimeder} (a) and \eqref{eq:cnvlem3est1} (a), we have\footnote{\ig By the density of the space $\C_{c}(\Rd;\R)$ in $\Ltp(\Rd;\R)$ we may assume $z\in\C_{c}(\Rd;\R)$.\fi  Clearly, $(\sKKS_{0} z)(t,x)$ may trivially be turned into a sequence by letting $t_{n}=1/n$, then $n\nearrow\infty$ means $t_{n}\searrow0$.}  
\ben\rencomrom
\item  the pointwise almost everywhere convergence $$(\sKKS_{0} z)(t,x)\underset{t\rightarrow0}{\longrightarrow}z(x)\mbox{ a.e., }$$ 
which follows by a standard argument since $\KKStxy$ is the fundamental solution of \eqref{eq:lkspde}; and 
 \item the convergence of the $\Ltp$ norms $$|(\sKKS_{0} z)(t))|_{2p}\underset{t\rightarrow0}{\longrightarrow}\left|z\right|_{2p},$$ 
 which follows from (i) together with the general Vitali convergence theorem (e.g., \cite[Chapter 5, p.~98]{RoyFitz}), since the family $\{|(\sKKS_{0} z)(t,x)|^{2p}\}_{t>0}$ can readily be checked to be both uniformly integrable and tight.
\een
The last two types of convergence in (i) and (ii) and an application of Fatou's lemma then give us the continuity at time $0$ in the $\Ltp$ norm (e.g., \cite[Theorem 7, p.~148]{RoyFitz}): 
\beq\lbl{eq:Vitali}
\bsp
\left|(\sKKS_{0} z)(t)-(\sKKS_{0} z)(0)\right|_{2p}=\left|(\sKKS_{0} z)(t)-z\right|_{2p}\longrightarrow 0 \mbox{ as }t\searrow0.
\end{split}
\eeq
The proof is complete.
\epf
\section{Local $\Ltp$ solution via the $N$-ball KSB approximations}\lbl{sec:local}
\subsection{Existence, uniqueness, and $\Ltp$ uniformity for KSB $N$-ball approximations}
Let $N\in\N$ be any fixed arbitrary natural number and, for $p\ge1$, let $d<6p$. Let $\B_{N}(0)\subset\Ltp(\Rd;\R)$ be the ball with radius $N$ centered at the origin in $\Ltp(\Rd;\R)$. Define the mapping $\sI_N^{(p)}:\Ltp(\Rd;\R)\mapsto \B_{N}(0)$ by
\beq\lbl{eq:locop}
\sI_N^{(p)}(v)= \begin{cases}
  v, & \text{if } |v|_{2p}\leq N, \vspace{2.5mm}\\ 
  \df{v}{|v|_{2p}}N, & \text{if }|v|_{2p}> N.
\end{cases}
\eeq

For each $N\in\N$, we first show the existence and uniqueness of a uniform $\Ltp$ solution for the following $\sI_N^{(p)}$-approximation of the Kuramoto-Sivashinsky Burgers equation and its L-KS formulation \eqref{eq:ibtbapsol}:
\beq\lbl{eq:sinks}
\bsp
U_{N}(t,x)&=\intrd\KKStxy \uny dy
+\tfrac1{2}\intrdzt\KKStsxy(\sI_N^{(p)}U_{N})(s,y)ds dy
\\&+\tfrac12\int_{\Rd}\int_{0}^{t}\sum_{\l=1}^{d}\frac{\pa\KKStsxy}{\pa y_{\l}} [(\sI_N^{(p)} U_{N})(s,y)]^2 dsdy.
\end{split}
\eeq

\bthm[Existence, uniqueness, and $\Ltp$ uniformity for the $\sI_N^{(p)}$-KS equation \eqref{eq:sinks}]\lbl{thm:unqexloc}  Suppose $p\ge1$ and $d<6p$.  Assume further that  $\un\in \Ltp\lpa\Rd;\R\rpa$.  Then, for any fixed $N\in\N$, there exists a unique solution $U_{N}\in\C\lpa[0,T];\Ltp\lpa\Rd;\R\rpa\rpa$ to \eqref{eq:sinks} such that
\beq\lbl{eq:unifbdloc}
\sup_{t\in [0,T]}|U_{N}(t)|_{2p}\leq C\lpa\lab u_{0}\rab_{2p}+TN+d N^2 T^{\frac{6p-{d}}{8p}}\rpa<\infty.
\eeq
\ethm
\thmref{thm:unqexloc} holds for every $N\in\N$ and will be used  to prove uniqueness and local existence for the KSB equation \eqref{eq:lkspdea} (equivalently \eqref{eq:ksBur}) in \thmref{thm:unqlocex2} below.  
The proof of \thmref{thm:unqexloc} follows from the next two lemmas.  

\blm[Uniform $\Ltp$ norm for the $\sI_N^{(p)}$-KS operator in \eqref{eq:sinks}]\lbl{lm:solunifl2bound}
Suppose $p\ge1$, $d<6p$, and $\un\in \Ltp\lpa\Rd;\R\rpa$. Fix $N\in\N$ and, for $U_{N}\in\MsB\lpa[0,T];\Ltp\lpa\Rd;\R\rpa\rpa$\footnote{The space of Borel-measurable, $\Ltp\lpa\Rd;\R\rpa$ valued, functions on $[0,T]$ (\appref{sec:acrnot} (10)).},  let the operator $\sO_{N}$ be given by
\beq\lbl{eq:Oop}
\lpa\mathscr O_{N} U_{N}\rpa(t,x)=(\sKKS_{0} \un)(t,x)+\sum_{k=1}^{2}\lpa{\mathscr  O_{k,N}} U_{N}\rpa(t,x),
\eeq
where 
\beqs
\bsp
(\sKKS_{0} \un)(t,x)&=\int_{\R^d}\KKStxy u_0(y) dy\\
(\mathscr O_{1,N} U_{N})(t,x)&=\tfrac12\lpa\sKKS (\sI_N^{(p)} U_{N})\rpa(t,x)
\\&=\tfrac1{2}\intrdzt\KKStsxy(\sI_N^{(p)}U_{N})(s,y)ds dy,
\\(\sO_{2,N}U_{N})(t,x)&=\tfrac12\sum_{\l=1}^{d}\lpa\sKKSl \lpa[\sI_N^{(p)} U_{N}]^{2}\rpa\rpa(t,x)
\\&=\tfrac12\sum_{\l=1}^{d}\int_{\Rd}\int_{0}^{t}\frac{\pa\KKStsxy}{\pa y_{\l}} [(\sI_N^{(p)} U_{N})(s,y)]^2 dsdy.
\end{split}
\eeqs
Then, 
\beq\lbl{eq:unballbound}
\sup\limits_{t\in [0,T]}|(\mathscr O_{N} U_{N})(t)|_{2p}\leq C\lpa\lab u_{0}\rab_{2p}+TN+d N^2 T^{\frac{6p-{d}}{8p}}\rpa<\infty,
\eeq  
for some constant $C\in(0,\infty)$.  Moreover,  $\mathscr O_{N} U_{N}\in\C\lpa[0,T];\Ltp\lpa\Rd;\R\rpa\rpa$.  
\elm
\bpf  First, by \lemref{lm:bdsKScnv} equation \eqref{eq:cnvlem2est1} (a) and by the definition of $\sI_N^{(p)}$ \eqref{eq:locop} we have 
\beq\lbl{eq:O1ub}
\bsp
\left|(\mathscr O_{1,N} U_{N})(t)\right|_{2p}&=\tfrac12\lab\lpa\sKKS (\sI_N^{(p)}U_{N})\rpa(t)\rab_{2p} \\&\le C\int\limits_0^t|(\sI_N^{(p)} U_{N})(r)|_{2p} dr\le C T N.
\end{split}
\eeq

Second, by Minkowski's inequality, \lemref{lm:Burgersfundest} \eqref{eq:cnvlem1est1} (a), and \eqref{eq:locop} we have
\beq\lbl{eq:O2ub}
\bsp
\left|(\mathscr O_{2,N} U_{N})(t)\right|_{2p}&\le\tfrac12\sum_{\l=1}^d \lab\lpa\sKKSl \lpa[\sI_N^{(p)} U_{N}]^{2}\rpa\rpa(t)\rab_{2p} 
\\&\leq C\sum_{\l=1}^d \int\limits_0^t (t-s)^{-\frac{2p+d}{8p}}\left|(\sI_N^{(p)} U_{N})(s)\right|_{2p}^{2} ds
\\&\leq C d N^2 T^{\frac{6p-{d}}{8p}}.
\end{split}
\eeq

Finally, by \eqref{eq:unzz} in the proof of \lemref{lm:bdsKSinitcnv},  \eqref{eq:cnvlem3est1} (a), we have
\beq\lbl{eq:un}
\bsp
\left|(\sKKS_{0} \un)(t)\right|_{2p} &\le\lab\KKStdot\rab_{1}\lab u_{0}\rab_{2p} \le C\lab u_{0}\rab_{2p} <\infty, \ 0<t\le T.
\end{split}
\eeq
By \eqref{eq:Oop}, \eqref{eq:O1ub}, \eqref{eq:O2ub}, \eqref{eq:un}, and Minkowski's inequality the desired uniform $\Ltp$ bound follows.  The conclusion $\mathscr O_{N} U_{N}\in\C\lpa[0,T];\Ltp\lpa\Rd;\R\rpa\rpa$ is a direct consequence of the convolution lemmas \ref{lm:Burgersfundest}--\ref{lm:bdsKSinitcnv} together with the definitions of $\mathscr O_{N} U_{N}$ and $\sI_N^{(p)} U_{N}$.
\epf
\blm[A fixed point lemma for the $\sI_N^{(p)}$-KS equation \eqref{eq:sinks}]\lbl{lm:fxdpt}
Suppose that the assumptions of \thmref{thm:unqexloc} are in force.   Let $\U_{\lambda}$ denote the Banach space of $\MsB\lpa[0,T];\Ltp\lpa\Rd;\R\rpa\rpa$ functions $U$ such that $u(0)=u_0$ equipped with the norm
\beq\lbl{eq:lambdanorm}
 \lab U\rab_{\U_{\lambda}}=\lbk\int\limits_0^T e^{-\lambda t}|U(t)|_{2p}^{2p}dt\rbk^{1/2p}<\infty, \ \lambda\in(0,\infty).
\eeq
Then, for each $N\in\N$, denoting by $\Gamma$ the standard Gamma function, there exists a constant $C\in(0,\infty)$ that may depend on $N$ and $T$ and a $\lambda_{0,d}>0$ satisfying
\beq\lbl{eq:lambda0}
\lbk{{\lambda_{0,d}^{-1}}+\Gamma \left( \tf{6p-d}{8p} \right) \lpa\lambda_{0,d}\rpa^{\tf{d-6p}{8p}}}\rbk=1\wedge\df1C, \ d<6p,
\eeq
 such that for any $\lambda>\lambda_{0,d}$
the operator $\mathscr O_{N}$ in \eqref{eq:Oop} is a contraction on $\U_{\lambda}$ and hence there exists a unique solution for \eqref{eq:sinks}, $U_{N}\in\U_{\lambda}$.  Moreover, $U_{N}\in\C\lpa[0,T];\Ltp\lpa\Rd;\R\rpa\rpa$ and $U_{N}$ satisfies the uniform $\Ltp$ norm bound in \eqref{eq:unifbdloc}.
\elm
\bpf
Fix $\lambda\in(0,\infty)$.  
For each $N\in\N$, define $\sO_{N}$ on $\U_{\lambda}$ by \eqref{eq:Oop}.  From \lemref{lm:solunifl2bound} above, it follows that the operator $\sO_{N}$  maps the Banach space $\U_{\lambda}$ into itself. Fix two elements $U,V\in\U_{\lambda}$ and fix $N\in\N$.  

Proceeding as in the proof of \lemref{lm:bdsKScnv}, \eqref{eq:cnvlem2est1} part (a), we apply Minkowski's and Young's inequalities followed by an application of \lemref{lm:basicfacts} \eqref{eq:lksL1bdskertimeder} (a) and Jensen's inequality we get
\beqs
\bsp
&|\lpa\sO_{1,N} U\rpa(t)-\lpa\sO_{1,N} V\rpa(t)|_{2p}^{2p}\\&=\df12\int_{\Rd}\lab\int_0^t\lpa{\KKStrdot}*\lbk(\sI_N^{(p)} U)(r,\cdot)-(\sI_N^{(p)}V)(r,\cdot)\rbk\rpa(x) dr\rab^{2p}dx
\\&\leq\df12\lbk  \int_0^t\left|\lab \KKStrdot\rab*\lab (\sI_N^{(p)} U)(r,\cdot)-(\sI_N^{(p)}V)(r,\cdot)\rab\right|_{2p}dr\rbk^{2p}
\\&\leq\df12\lbk \int_0^t \left| \KKStrdot\right|_1 \left|(\sI_N^{(p)} U)(r,\cdot)-(\sI_N^{(p)}V)(r,\cdot)\right|_{2p} dr\rbk^{2p}
\\&\le C{T^{2p-1}} \int_0^t \left|(\sI_N^{(p)} U)(r,\cdot)-(\sI_N^{(p)}V)(r,\cdot)\right|_{2p}^{2p}  dr, 
\end{split}
\eeqs 
which, upon using \eqref{eq:locop}, implies the existence of a constant $C\in(0,\infty)$ (depending on $N$, $T$, and $p$) such that
\beq\lbl{eq:o1estfplm}
|\lpa\sO_{1,N} U\rpa(t)-\lpa\sO_{1,N} V\rpa(t)|_{2p}^{2p}
\leq C \int_0^t \left|U(s) - V (s) \right|_{2p}^{2p} ds. 
\eeq

Using the basic inequality $|\sum_{i=1}^{d}\lab a_{i}\rab|^{2p}\le C\sum_{i=1}^{d}| a_{i}|^{2p}$, for all $p\ge1/2$, \lemref{lm:Burgersfundest}, Jensen's inequality with respect to the probability measure on $[0,t]$ given by 
${\mu_{t}^{d}(ds)}/{|\mu_{t}^{d}|},$ where $\mu_{t}^{d}(ds)={(t-s)^{-\frac{2p+d}{8p}}ds}$ and $|\mu_{t}^{d}|=t^{\frac{6p-d}{8p}}$,
 H\"older's inequality, and the fact that $|\sI_N^{(p)}(v)|_{2p}\le N$, we get
\beqs
\bsp
&|\lpa\sO_{2,N} U\rpa(t)-\lpa\sO_{2,N} V\rpa(t)|_{2p}^{2p}
\\&\leq C\sum_{\l=1}^d \left|\int_0^t \int_{\Rd} \frac{\pa\KKtsSy}{\pa y_{l}}\lab(\sI_N^{(p)}  U)^{2}(s,y)-(\sI_N^{(p)}  V)^{2}(s,y)\rab dy ds\right|_{2p}^{2p}
\\&\leq C \lbk \int_0^t (t-s)^{-\frac{2p+d}{8p}}\lpa\int_{\Rd}\lab(\sI_N^{(p)}  U)^{2}(s,y)-(\sI_N^{(p)}  V)^{2}(s,y)\rab^{p} dy\rpa^{1/p}  ds\rbk^{2p}
\\&\leq C|\mu_{t}^{d}| ^{2p-1} \int_0^t \left|\lab(\sI_N^{(p)} U)(s)-(\sI_N^{(p)} V)(s)\rab[|(\sI_N^{(p)} U)(s)|+|(\sI_N^{(p)} V)(s)|]\right|_p^{2p} {\mu_{t}^{d}(ds)}
\\&\leq C_{N,T,p} \int_0^t (t-s)^{-\frac{2p+d}{8p}}\left|(\sI_N^{(p)} U)(s)-(\sI_N^{(p)} V)(s)\right|_{2p}^{2p} ds,
\end{split}
\eeqs
where we used the fact that $|\mu_{t}^{d}| \le T^{\frac{6p-d}{8p}}$.  This and the definition of $\sI_N^{(p)}$ in \eqref{eq:locop} yield
\beq\lbl{eq:o2estfplm}
|\sO_{2,N} U(t)-\sO_{2,N} V(t)|_{2p}^{2p}\leq C \int_0^t (t-s)^{-\frac{2p+d}{8p}} \left|U(s) - V (s) \right|_{2p}^{2p} ds. 
\eeq
From estimates \eqref{eq:o1estfplm} and \eqref{eq:o2estfplm} we deduce that
\beq\lbl{eq:contrest}
\bsp
&|\sO_{N} U-\sO_{N} V|_{\U_{\lambda}}^{2p}=\int_0^T e^{-\lambda t}|\sO_{N} U(t)-\sO_{N} V(t)|_{2p}^{2p} dt
\\&\leq C\int_0^T e^{-\lambda t}\lpa \int_0^t \lbk1+(t-s)^{-\frac{2p+d}{8p}}\rbk\left|U(s) - V (s) \right|_{2p}^{2p} ds\rpa dt
\\&\leq C\lpa \int_0^\infty e^{-\lambda \rho}\lbk1+\rho^{-\frac{2p+d}{8p}}\rbk d\rho \rpa |U-V|_{\U_{\lambda}}^{2p}
\\&=C\lbk{{\lambda}^{-1}+\Gamma \left( \tf{6p-d}{8p} \right) \lambda^{\tf{d-6p}{8p}}}\rbk|U-V|_{\U_{\lambda}}^{2p}.
\end{split}
\eeq
Picking the constant $C$ at the end of \eqref{eq:contrest}, let $\lambda_{0,d}$ be such that 
\beqs
\lbk{{\lambda_{0,d}^{-1}}+\Gamma \left(\tf{6p-d}{8p} \right) \lpa\lambda_{0,d}\rpa^{\tf{d-6p}{8p}}}\rbk=1\wedge\df1C, \ d<6p.
\eeqs
Then, for any  $\lambda>\lambda_{0,d}$, we have that the operator $\sO_{N}$ is a contraction on $\U_{\lambda}$. Consequently, for each $\lambda>\lambda_{0,d}$, there exists a unique fixed point $U_{N}=\sO_{N}U_{N}$ for this operator; i.e., there exists a unique solution $U_{N}\in\U_{\lambda}$ for \eqref{eq:sinks}.  By \lemref{lm:solunifl2bound},  $U_{N}\in\C\lpa[0,T];\Ltp\lpa\Rd;\R\rpa\rpa$ and satisfies the uniform $\Ltp$ norm bound in \eqref{eq:unifbdloc}. 
\epf

\bpfs{Completing the proof of \thmref{thm:unqexloc}} From \lemref{lm:solunifl2bound} and \lemref{lm:fxdpt} we deduce the desired existence, uniqueness, continuity, and uniform $\Ltp$ bound in \eqref{eq:unifbdloc} of the solution to \eqref{eq:sinks} as in \thmref{thm:unqexloc}.  
\epfs
\subsection{Global uniqueness and local existence for the KSB equation}
\thmref{thm:unqexloc} implies the uniqueness and local-in-time existence of $\Ltp$, $p\ge1$, solutions for the KSB PDE \eqref{eq:ksBur} in spatial dimensions $d<6p$.  This is the content of the next result, which restates our main result in \thmref{thm:exunreg}.
\bthm[Uniqueness and local existence for the KS Burgers PDE \eqref{eq:ksBur}]\lbl{thm:unqlocex2}  Suppose $p\ge1$ and $d<6p$.  Assume further that  $\un\in \Ltp\lpa\Rd;\R\rpa$.   Let $U_{1},U_{2}\in\C\lpa[0,T];\Ltp\lpa\Rd;\R\rpa\rpa$ be solutions to \eqref{eq:ksBur}.  Then, $U_{1}(t)=U_{2}(t)$  for all $t\in[0,T].$    Furthermore, there exists a solution $U\in\C\lpa[0,\tau_{\infty});\Ltp\lpa\Rd;\R\rpa\rpa$ to the KS PDE \eqref{eq:ksBur}, for some $0<\tau_{\infty}\le T$. 
\ethm
Our proof below borrows from, and has the flavor of, stopping-times arguments from stochastic analysis (e.g., \cite[equations (3.15) and (2.3), respectively]{Acom,Acom1} and \cite{DapZ}), adapted to the deterministic setting here.
\bpf  Let $U_{1}$ and $U_{2}$ be as in \thmref{thm:unqlocex2} above, then they satisfy the KS \eqref{eq:ibtbapsol}.  Now, define the  times 
$$\sigma_{N}:= T\wedge\inf\lbr t\ge0;\lab U_{1}(t)\rab_{2p}\wedge\lab U_{2}(t)\rab_{2p}\ge N\rbr,\quad N\in\N;$$ 
then the functions
\beqs
U_{i}^{(N)}=\lbr U_{i}(t\wedge\sigma_{N});0\le t\le T\rbr,\ i=1,2,
\eeqs
are solutions to the approximating $\sI_N^{(p)}$-KS integral equation \eqref{eq:sinks}, for each $N\in\N$.  So, \thmref{thm:unqexloc} implies that 
\beqs
\bsp
U_{1}(t\wedge\sigma_{N})&=U_{2}(t\wedge\sigma_{N}), \ \forall (N,t)\in\N\times[0,T].
\end{split}
\eeqs
Taking the limit as $N\nearrow\infty$ and using the time continuity of $U$ and $V$ and the $\sigma_{N}$ definition, we get 
\beqs
\bsp
U_{1}(t)=\lim_{N\nearrow\infty}U_{1}(t\wedge\sigma_{N})=\lim_{N\nearrow\infty}U_{2}(t\wedge\sigma_{N})=U_{2}(t),\ t\in[0,T]; 
\end{split}
\eeqs 
and uniqueness is established for \eqref{eq:ksBur}.

For the local existence assertion, let $U_{N}$ be the solution to the  $\sI_N^{(p)}$-KS approximation \eqref{eq:sinks} for every $N\in\N$.  Define the time $\tau_{N}$ by
$$\tau_{N}:=T\wedge\inf\lbr t\ge0;\lab U_{N}(t)\rab_{2p}\ge N\rbr,\quad N\in\N.$$ 
Clearly, 
\beq\lbl{eq:consis}
U_{M}(t)=U_{N}(t)\mbox{ for all $M\ge N$ and $t\le \tau_{N}$}. 
\eeq 
Moreover, the function $\lbr U_{N}(t); 0\le t\le \tau_{N} \rbr$
satisfies the mild L-KS formulation \eqref{eq:ibtbapsol}, and hence the KSB PDE, on the time interval $[0,\tau_{N}]$.  Let $\tau_{\infty}=\sup_{N\nearrow\infty}\tau_{N}\le T$.  Then, using the consistency property in \eqref{eq:consis} and setting $U(t):=U_{N}(t)$ for $t\le\tau_{N}$ and $N\ge1$ we obtain a solution $U$ to the KS PDE \eqref{eq:ksBur} on the time interval $[0,\tau_{\infty})$\footnote{As we noted in \secref{sec:mainsketch}, we may think of the just-constructed $U$ as the function whose path is obtained by glueing the path of $U_{N}$ onto that of $U_{N-1}$ for all $N\ge2$.}.   
\epf
\brm[Added Local existence and uniqueness for the KS PDE variant]\lbl{rem:locexunqmks} The same statements and conclusions of \thmref{thm:unqexloc} and \thmref{thm:unqlocex2} hold unchanged for the $\sI_N^{(p)}$ version of the KS variant \eqref{eq:kspdeavar}:
\beq\lbl{eq:locksvar}
\bsp
U_{N}(t,x)&=\intrd\KKStxy \uny dy
\\&+\tfrac12\int_{\Rd}\int_{0}^{t}\sum_{\l=1}^{d}\frac{\pa\KKStsxy}{\pa y_{\l}} [(\sI_N^{(p)} U_{N})(s,y)]^2 dsdy
\end{split}
\eeq
and for \eqref{eq:kspdeavar} itself, respectively.  The proofs are the same as above except that the extra factor 
$$\tfrac1{2}\intrdzt\KKStsxy(\sI_N^{(p)}U_{N})(s,y)ds dy$$ 
along with the need for \lemref{lm:bdsKScnv} disappear.
\erm

\appendix
\section{Briefly on Brownian-time Brownian motion connection\lbl{app:briefdisc}}
As we remarked in our earlier work \cite{Alksspde,Abtpspde,Aks,AX17}, there is an intimate connection between the L-KS kernel and the Brownian-time Brownian motion (BTBM) kernel:
\beq\lbl{eq:btbmk}
\KBtxy=2\int_{0}^{\infty}\psxy\ptsz ds=2\int_{0}^{\infty}\frac{\e^{-|x-y|^{2}/2s}\e^{-s^{2}/2t}}{\lpa2\pi s\rpa^{d/2}\sqrt{2\pi t}}ds,
\eeq
which is both
\ben\rencomrom
\item the density of a BTBM $\{X^{x}(|B_{t}|);t\ge0\}$, where $X^{x}$ is a $d$-dimensional Brownian motion starting at $x\in\Rd$ and $B$ is an independent one-dimensional BM starting at $0$; and 
\item the fundamental solution to both the time-fractional and the fourth order PDE pair: 
\beq\lbl{eq:BTBMhalf}
\bsp
\pa_{t}^{1/2}u=\tfrac1{\sqrt 8}\lap u\  \mbox{ and }\  \pa_{t}u=\tfrac{1}{\sqrt{8\pi t}}\lap \un+\tf18\lap^{2}u,\mbox{ with }
u(0,x)=\delta(x),
\end{split}
\eeq
which are equivalent for suitable initial conditions $u(0,x)$ (see \cite{Atfhosie,AZ01} and the references therein for more details).
\een
  The cozy relationship between the BTBM and L-KS kernels, and their PDEs, is due to the underlying Brownian-time construction in both (see \cite[equation (2.12) and Lemma 3.1]{Alksspde} for similar estimates for the two kernels).  Using standard heat kernel estimates and the definition of $\KBtx$ in 
\eqref{eq:btbmk}, we can easily prove that \lemref{lm:basicfacts} \eqref{eq:lksfundbdskertimeder} and hence \eqref{eq:lksL1bdskertimeder} both hold for $\KBtx$ in place of the L-KS kernel $\KKStx$, with different constants\footnote{Of course, for the BTBM density $|\KBtdot|_{1}=1$.}.  
\blm[Twin BTBM kernel fundamental estimates]\lbl{lm:btbmlksfest}
 The BTBM density \eqref{eq:btbmk} satisfies the same estimates as those for the L-KS kernel in \eqref{eq:lksfundbdskertimeder}, with only different constants.  I.e.,
\beq\lbl{eq:btbmbds}
\bsp
&(a) \lab\KBtxy\rab= 2\int_{0}^{\infty}\frac{\e^{-|x-y|^{2}/2s}\e^{-s^{2}/2t}}{\lpa2\pi s\rpa^{d/2}\sqrt{2\pi t}}ds\\
&(b) \lab\frac{\partial\KBtxy}{\partial y_\l}\rab\le C_{2}\int_{0}^{\infty}\frac{e^{-{c_{2}|x-y|^2}/{s}}}{s^{\tf{
d+1}{2}}} \frac{e^{-\frac{s^2}{2 t}}}{\sqrt{2\pi t}} d{s}
\\&(c) \lab\frac{\partial\KBtxy}{\partial t}\rab\le C_{3} \int_{0}^{\infty}\frac{e^{-{c_{1}|x-y|^2}/{s}}}{s^{{d}/{2}}}\frac{e^{-\frac{c_{3}s^2}{t}}}{t^{3/2}} ds
\\&(d) \lab\frac{\partial^{2}\KBtxy}{\partial t\partial y_\l}\rab\le C_{4}\int_{0}^{\infty}\frac{e^{-\frac{c_{2}|y|^2}{s}}}{s^{({d+1})/{2}}}\frac{e^{-\frac{c_{3}s^2}{t}}}{t^{3/2}} ds, 
\end{split}
\eeq
for  $(t,x,y)\in(0,T]\times\Rd\times\Rd$ and some constants $C_{i},c_{i}>0$, $i=1,2,3$ and $C_{4}=C_{2}C_{3}$, that may only depend on $T$ and $d$.
\elm
\bpf
Part (a) is is just \eqref{eq:btbmk}.  Part (b) follows easily from \eqref{eq:btbmk} and the following standard estimate for the $d$-dimensional Brownian motion density
\beq\lbl{eq:bmspder}
\left|\frac{\partial}{\partial y_\l}\psxy\right|\leq C_{2}\frac{e^{-\frac{c_{2}|x-y|^2}{s}}}{s^{({d+1})/{2}}};\   s\in(0,\infty), x,y\in\Rd, \mbox{ and }\l=1,\ldots,d,
\eeq
while part (c) follows readily from \eqref{eq:btbmk} and the following standard time-derivative estimates for the one dimensional BM density
\beq\lbl{eq:bmtimeder}
\left|\frac{\partial}{\partial t}\ptsz\right|\leq C_{3}\frac{e^{-\frac{c_{3}s^2}{t}}}{t^{3/2}}; \ 0 < t < T \mbox{ and } y\in\Rd.
\eeq 
The estimate for the mixed spatio-temporal derivatives now   immediately follows.  Denote by $\pa_{\l}$ and $\pa_{t}$ the partial derivatives $\pa/\pa y_{\l}$ and $\pa/\pa{t}$, respectively.  We then easily have from \eqref{eq:btbmk}, \eqref{eq:bmspder}, and \eqref{eq:bmtimeder} that
\beq\lbl{eq:kermixedder}
\bsp
\lab\frac{\partial^{2}\KBtxy}{\partial t\partial y_\l}\rab&\le2\int_{0}^{\infty}\lab\pa_{\l}\psxy\rab\lab\pa_{t}\ptsz\rab ds
\\&\le C\int_{0}^{\infty}\frac{e^{-\frac{c_{2}|x-y|^2}{s}}}{s^{({d+1})/{2}}}\frac{e^{-\frac{c_{3}s^2}{t}}}{t^{3/2}} ds,
\end{split}
\eeq
completing the proof.
\epf
\section{Briefly on weak stability and space-time continuity}\lbl{app:spatialcont}
We start by rewriting the L-KS mild formulation \eqref{eq:ibtbapsol} of the KSB PDE in terms of the operators $\sKKS_{0}$, $\sKKS$, and $\sKKSl$ in \eqref{eq:convlks}: 
\beq\lbl{eq:LKSoperform}
\bsp
U(t,x)&=(\sKKS_{0} \un)(t,x)+\tf12(\sKKS U)(t,x)+\tf12\ds\sum_{\l=1}^{d}
(\sKKSl U^{2})(t,x) 
\end{split}
\eeq

\subsection{On weak $\Ltp$ stability}
\bthm[Global weak $\Ltp$ stability for KSB]\lbl{thm:wstab}
Assume that $T>0$ and $p\ge1$ are arbitrary and fixed, and let $d<6p$.   If $\uno,\unt\in \Ltp\lpa\Rd;\R\rpa$ and if $\Uk\in\C\lpa[0,T];\Ltp\lpa\Rd;\R\rpa\rpa$ is the unique L-KS solution to the KSB PDE \eqref{eq:ksBur} corresponding to $
\unk$,  $k=1,2$, then
\beq\lbl{eq:wstabinq}
\lab\Uo(t)-\Ut(t)\rab_{2p}\le C_{T,d,p}\lab\uno-\unt\rab_{2p},
\eeq
where the constant $C_{T,d,p}\in(0,\infty)$ depends on $T,d,p$ as well as on $$\max_{k=1,2}\sup_{0\le r\le T}\lab\Uk(r)\rab_{2p}.$$
In particular, the weak stability in \eqref{eq:wstabinq} holds for the KSB PDE for $t\in[0,\tau]$, where $\tau$ is as in \thmref{thm:exunreg}.
\ethm
\brm\lbl{rem:wstab}   We remind the reader that our result in \thmref{thm:exunreg} asserts  global uniqueness, so we assume only the existence of global solutions $\Uk$, $k=1,2$, in \thmref{thm:wstab}. The dependence on the norms 
$$\max_{k=1,2}\sup_{0\le r\le T}\lab\Uk(r)\rab_{2p}$$
is why we call the stability weak.
\erm 
\bpf
Using \eqref{eq:LKSoperform} together with \lemref{lm:Burgersfundest}--\lemref{lm:bdsKSinitcnv}, and then Cauchy-Schwarz inequality we readily get
\beq\lbl{eq:stabderivation}
\bsp
&\lab\Uo(t)-\Ut(t)\rab_{2p}\le\lab\lpa\sKKS_{0} \lbk\uno-\unt\rbk\rpa(t)\rab_{2p}
\\&\quad +\tf12\lab\lpa\sKKS \lbk\Uo-\Ut\rbk\rpa(t)\rab_{2p}
\\&\quad +\tf12\ds\sum_{\l=1}^{d}\lab\lpa\sKKSl \lbk\lpa\Uo\rpa^{2}-\lpa\Ut\rpa^{2}\rbk\rpa(t)\rab_{2p}
\\&\le C\lbr\lab\uno-\unt\rab_{2p}+\int_{0}^{t}\lab\Uo(r)-\Ut(r)\rab_{2p}dr\rbr
\\&\quad +C\int_{0}^{t}(t-r)^{-\tf{2p+d}{8p}}\lab\lbk\Uo(r)-\Ut(r)\rbk \lbk\Uo(r)+\Ut(r)\rbk\rab_{p}dr
\\&\le C\int_{0}^{t}\lab\Uo(r)-\Ut(r)\rab_{2p} \lbk1+(t-r)^{-\tf{2p+d}{8p}}\max_{k=1,2}\sup_{0\le r\le T}\lab\Uk(r)\rab_{2p}\rbk dr
\\&\quad+C\lab\uno-\unt\rab_{2p}.
\end{split}
\eeq
The desired conclusion now follows from Gr\"onwall's inequality.
\epf

\subsection{On continuity in space-time}
The continuity in \remref{rem:spatialcont} follows easily from that of each term in \eqref{eq:LKSoperform}, which are clearly continuous for step functions.  We now  state and show a slightly stronger result for the operator $\sKKSl$, corresponding to the Burgers term, and we leave the simpler cases of $\sKKS$ and $\sKKS_{0}$ to the interested reader.
\blm\lbl{lem:space-time cont}
Fix arbitrary $T>0$ and $p>1$, and let $d<3p$.  The operator $\sKKSl$ is bounded from $\Lgam\lpa[0,T];\Lp\lpa\Rd\rpa\rpa$ into $\C_{b}\lpa[0,T]\times\Rd;\R\rpa$, and it satisfies the following estimate
\beq\lbl{eq:cnvlem1est2}
\lab\sKKSl(u)(t,x)\rab\le\int_{0}^{t}\lpa t-r\rpa^{-\frac {p+d}{4p}}|u(r,\cdot)|_{p} dr\le C\left(\int_0^t|u(r)|_p^\gamma dr\right)^{1/\gamma},
\eeq
for $\gamma>4p/(3p-d)$.
\elm
\bpfs{Proof}
 Two H\"older inequality applications in space then in time;  and the estimate in \eqref{eq:lqderkergen}  (see also \eqref{eq:Parskerderl2}, \eqref{eq:kerderl2q1}, \eqref{eq:kerderl2q2}) with $q=p/(p-1)$, give\footnote{We note here that in \eqref{eq:lqderkergen} the $(p,q)$ where a Young-pair not a H\"older-pair like in \eqref{eq:spctimeHol} above.} 
\beq\lbl{eq:spctimeHol}
\bsp
\left|\sKKSl(u)(t,x)\right|&\le\int_0^t \left|\partial_\l \KKStrxdot\right|_q|u(r)|_{p} dr
\\&\le C\int_{0}^{t}\lpa t-r\rpa^{-\frac {p+d}{4p}}|u(r)|_{p} dr
\leq C\left(\int_0^t|u(r)|_p^\gamma dr\right)^{1/\gamma},
\end{split}
\eeq
for $\gamma>4p/(3p-d)$, which establishes  \eqref{eq:cnvlem1est2}.  Clearly, $\{\sKKSl(u)(t,x); (t,x)\in[0,T]\times\Rd\}$ is continuous for step functions $u$, which completes the proof.  
\epfs
\section{Frequent acronyms and notations key}\lbl{sec:acrnot}
\begin{enumerate}\renewcommand{\labelenumi}{\Roman{enumi}.}
\item {\textbf{Acronyms}}\vspace{2mm}
\begin{enumerate}\renewcommand{\labelenumii}{(\arabic{enumii})}
\item BM: Brownian motion.
\item BTBM: Brownian-time Brownian motion.
\item KS: Kuramoto-Sivashinsky.
\item KSB: Kuramoto-Sivashinsky-Burgers.
\end{enumerate}
\vspace{2.5mm}
\item {\textbf{Notations}}\vspace{2mm}\\
Throughout the article sets of numbers or Euclidean spaces are denoted using the {\tt{$\backslash$mathbb}} font  (e.g., $\R$, $\Rd$, $\Cmpl$, $\N$, etc...); while spaces of functions are denoted using the {\tt{$\backslash$mathds}} font (e.g., $\C(\Rd;\R)$, $\Lp(\Rd;\R)$, etc...).  The following is a list of some notation we use in this article:
\begin{enumerate}\renewcommand{\labelenumii}{(\arabic{enumii})}
\item $\KKStxy$ is the L-KS kernel;
\item For suitable functions $z,u$, and $v$ (see \eqref{eq:convlks}) the operators $\sKKS_{0}$, $\sKKS$, and $\sKKSl$ are given by:
\beqs
\bsp
(\sKKS_{0} z)(t,x)&= \int_{\R^d}\KKStxy z(y)dy,\\
(\sKKS u)(t,x)&= \int_0^t\int_{\R^d}\KKStrxy u(r,y)dy dr,\\
(\sKKSl v)(t,x)&=\int_0^t\int_{\R^d}\frac{\partial \KKStrxy}{\partial y_\l}v(r,y)dy dr,\ \l=1,\ldots, d;
\end{split}
\eeqs
\item $*$ denotes the standard convolution operator;
\item $\FKKStxi$ is the Fourier transform of the L-KS kernel;  
\item $\hat u(\xi)=\lpa\sF u\rpa(\xi)$ is the Fourier transform of the function $u$, at $\xi\in\Rd$;
\item $\Lp(\Rd;\R):=\lbr f:\Rd\to\R;|f|_{p}^{p}=\int_{\Rd}\lab f(x)\rab^{p}dx<\infty\rbr$;
\item $|f|_{p}$ and $|f(\cdot)|_{p}$ denote the $\Lp(\Rd;\R)$-norm of a real function $f$;
\item $\C\lpa[0,T];\Lp\lpa\Rd;\R\rpa\rpa$ denotes the set of continuous functions $\{u(t);0\le t\le T\}$ that are $\Lp\lpa\Rd;\R\rpa$-valued; 
\item $\C_{b}\lpa[0,T]\times\Rd;\R\rpa$ denotes the set of real-valued continuous functions on $[0,T]\times\Rd$ that are bounded;
\item $\MsB\lpa[0,T];\Lp\lpa\Rd;\R\rpa\rpa$ denotes the set of Borel-measurable functions $u=\{u(t);0\le t\le T\}$ that are $\Lp\lpa\Rd;\R\rpa$-valued; 
\item $\C_{c}(\Rd;\R)$ is the set of continuus functions with compact support;
\item $\B_{N}(0)\subset\Ltp(\Rd;\R)$ is the ball with radius $N$ centered at the origin in $\Ltp(\Rd;\R)$.
\item $\pa_{\l}$ is the partial derivative in the $\l$-th spatial variable, $\l\in\{1,\ldots,d\}$;
\item $\pa_{t}$ is the partial derivative in the time variable $t$;
\item $\N=\{1,2,3,\ldots\}$ and $\N_{0}=\{0,1,2,3,\ldots\}$; and
\end{enumerate}
\end{enumerate}

\end{document}